\documentclass{amsart}
\usepackage{amssymb,color,mathrsfs, manfnt, wasysym,cancel
}

    \def\cD{\mathcal D} \def\cE{\mathcal E}    \def\cI{\mathcal I}           \def\cS{\mathcal S}  \def\cT{\mathcal T}  

\def\sF{\mathscr F}

\def\fx{\mathfrak{x}}\def\fy{\mathfrak{y}}\def\fz{\mathfrak{z}}\def\fZ{\mathfrak{Z}}
\def\fX{\mathfrak{X}}\def\fY{\mathfrak{Y}}\def\fm{\mathfrak{m}}\def\fp{\mathfrak{p}}\def\fq{\mathfrak{q}}
\def\fd{\mathfrak{d}}\def\fa{\mathfrak{a}} \def\fb{\mathfrak{b}}\def\fc{\mathfrak{c}}
\def\fp{\mathfrak{p}}\def\fA{\mathfrak{A}}\def\fB{\mathfrak{B}}\def\fC{\mathfrak{C}}
  
\def\fn{\mathfrak{n}} 
\def\Z{{\mathbb Z}}   \def\N{{\mathrm N}} \def\Q{{\mathbb Q}}

\def \ring {\Z[\omega]}

\DeclareMathOperator{\li}{Li}
\DeclareMathOperator{\tr}{Tr}

\DeclareMathOperator{\res}{Res}

\DeclareMathOperator{\re}{Re}

\def \e{e}
\def\cond#1{n_{#1}}

\renewcommand{\mod}[1]{\,{\rm mod}\,#1}

\def\su#1{\sum_{\substack{#1}}}
\def\ps#1{\sideset{}{^*}\sum_{\substack{#1}}}
\def\pr#1{\prod_{\substack{#1}}}

\def\bs#1{\begin{equation*} \begin{split} #1 \end{split} \end{equation*}}
\def\bsc#1{\begin{equation} \begin{split} #1 \end{split} \end{equation}}
\def\eqs#1{\begin{equation*} #1 \end{equation*}}
\def\eqn#1{\begin{equation} #1 \end{equation}}
\def\mult#1{\begin{multline*}#1\end{multline*}}
\def\multn#1{\begin{multline}#1\end{multline}}
\def\ar#1{\left\{ \begin{array}{l@{\quad\text{if }}l} #1 \end{array}\right.
}

\def\({\left(} \def\){\right)} \def\[{\left[} \def\]{\right]} 
\def\fl#1{\left\lfloor#1\right\rfloor} \def\ceil#1{\left\lceil#1\right\rceil}
\def\le{\leqslant} \def\ge{\geqslant}
\def\eps{{\varepsilon}}

\definecolor{pink}{rgb}{1,.2,.6}
\definecolor{orange}{rgb}{0.7,0.3,0}
\definecolor{blue}{rgb}{.2,.6,.75}
\definecolor{green}{rgb}{.4,.7,.4}
\definecolor{mycolor}{RGB}{210,15,35}
\def\red#1{\textcolor{mycolor}{#1}}

\newcommand*\cube{\mbox{\mancube}}

\newtheorem{lem}{Lemma}[section]
\newtheorem{prop}[lem]{Proposition}
\newtheorem{thm}[lem]{Theorem}
\newtheorem{defn}[lem]{Definition}

\newtheorem{theorem}[lem]{Theorem}
\theoremstyle{definition}

\newtheorem{rem}[lem]{Remark}
\newtheorem{remark}[lem]{Remark}

\def \w {5}

\begin{document}
\author{Ahmet M. G\"ulo\u{g}lu  }
\address{Department of Mathematics, Bilkent University, Ankara, Turkey}
\email{guloglua@fen.bilkent.edu.tr}
\keywords{Non-vanishing, Mollified Moments, Cubic Dirichlet characters, Cubic Gauss sums, Hecke $L$-functions.}

\author{Hamza Yesilyurt}
\address{Department of Mathematics, Bilkent University, Ankara, Turkey}
\email{hamza@fen.bilkent.edu.tr}

\thanks{The first author is supported by T\"UB\.ITAK Research Grant no. 119F413}

\title[Mollified Moments of Cubic Dirichlet L-Functions]{Mollified Moments of Cubic Dirichlet L-Functions over the Eisenstein Field}

\begin{abstract}
We prove, assuming the generalized Riemann Hypothesis (GRH) that there is a positive density of $L$-functions associated with primitive cubic Dirichlet characters over the Eisenstein field that do not vanish at the central point $s=1/2$. This is achieved by computing the first mollified moment, which is obtained unconditionally, and finding a sharp upper bound for the higher mollified moments for these $L$-functions, under GRH. The proportion of non-vanishing is explicit, but extremely small. 
\end{abstract}

\maketitle

\section{Introduction}
S. Chowla conjectured in \cite{Chow} that, for any real non-principal Dirichlet character $\chi, L(1/2, \chi) \neq 0$. Li \cite{Li} showed that over the rational function ﬁeld, there are infinitively many quadratic Dirichlet $L$-functions such that $L(1/2, \chi) = 0$. It is believed, however, that the number of such characters should be of density zero among all quadratic characters. 
For quadratic Dirichlet $L$-functions, Özlük and Snyder \cite{OS-1999} showed, under the Generalized Riemann Hypothesis (GRH), that at least 15/16 of the $L$-functions attached to quadratic characters do not vanish, by computing the one-level density for the low-lying zeroes in the family. The conjectures of Katz and Sarnak \cite{KS-book} imply that $L(1/2, \chi)  \neq 0$ for almost all quadratic Dirichlet $L$-functions. By computing the first two mollified moments Soundararajan \cite{Sound} proved without assuming GRH that at least 87.5\% of the quadratic Dirichlet $L$-functions do not vanish at $s = 1/2$. It is well-known that using the first two (non-mollified) moments does not lead to a positive proportion of non-vanishing, as they grow too fast (see Conjecture 1.5.3 in \cite{CFKRS} and the work of Jutila \cite{Jut}). In the function field case, Bui and Florea \cite{BuiFlo2018} obtained a proportion of non-vanishing of at least 94\% for quadratic Dirichlet $L$-functions by computing the one-level density.

In this paper, we consider the family of $L$-functions attached to primitive cubic Dirichlet characters defined over the Eisenstein field $K = \Q(\omega)$, where $\omega = e^{2\pi i/3}$. Many new
conceptual and technical difficulties appear when considering cubic characters, and the results in the literature are fewer and weaker. Several authors studied the one level density for families of cubic characters over both $\Q$ and $K$ for test functions whose Fourier transforms are supported in $(-1, 1)$. These families have unitary symmetry. For a family with unitary symmetry, one needs to go beyond $(-1,1)$ to obtain a positive proportion of non-vanishing. This was achieved for a thin subfamily of cubic characters over $K$ by David and the first author \cite{DG}, where the test functions are supported in $(-13/11, 13/11)$, which in turn gives $2/13$ as the proportion of non-vanishing. These results are obtained under GRH. 
As for the moments, the first moment of $L( 1/2, \chi)$, where $\chi$ is a primitive cubic character, was computed by Baier and Young \cite{BaYo} over $\Q$, by Luo \cite{Luo} for the same thin sub-family as in \cite{DG} over $K$, and by David, Florea, and Lalin \cite{DFL-1} over function fields. In these three papers, the authors obtained lower bounds for the number of non-vanishing cubic twists, but not positive proportions, by using upper bounds on higher moments. The lower bound in \cite{DFL-1} was later improved by Ellenberg, Li, and Shusterman in \cite{ELS2020}. Obtaining an asymptotic for the second moment for cubic Dirichlet $L$-functions is still an open question over functions fields or number fields. Recently, David, Florea and Lalin \cite{DFL-2} proved using mollified moments of $L$-functions that there is a positive proportion of non-vanishing for cubic Dirichlet $L$-functions at $s = 1/2$ over function fields in the non-Kummer case. This is obtained by computing the ﬁrst molliﬁed moment, and by finding a sharp upper bound for the second molliﬁed moment. Using the same approach, we prove the following result. 

\begin{thm} \label{thm:posprop}
Let $\sF$ be the family of primitive cubic Dirichlet characters defined in \eqref{family}. Then, assuming GRH for the $L$-functions in $\sF$, 
\eqs{\lim_{X \to \infty} \frac{\# \{ \chi \in \sF : L(1/2,\chi ) \neq 0, \quad n_\chi \le X \}}{\# \{ \chi \in \sF : n_\chi \le X \} } \ge \frac{3e^{-e^{101.3}}}{(\sqrt 3 - 1)^2}, }
where $n_\chi$ denotes the norm of the conductor of $\chi$.
\end{thm} 
As in \cite{DFL-2}, our positive proportion of non-vanishing is extremely small, since the implied constant in the upper bound for the second moment is too large. 

\subsection{Statement of Results}

In what follows, we shall write $\sF(X)$ to denote characters $\chi \in \sF$ with $n_\chi \le X$. Our first result gives an upper bound for the $k$\textsuperscript{th} mollified moment of the $L$-functions for any real $k > 0$. 

\begin{theorem}\label{thm:upperbound}
Let $k, \kappa > 0$ be real numbers such that $k\kappa$ is an integer bounded by an absolute constant. Let $\sF$ be the family of primitive cubic Dirichlet characters defined in \eqref{family}. Then, assuming GRH for the $L$-functions in $\sF$, the inequality 
\eqn{
\label{momentupperbound}
\sum_{\chi \in \sF(X)} |L(1/2,\chi)|^k |M(\chi, \kappa)|^{k\kappa} \ll_{k,\kappa} |\sF(X)|}
holds for $X$ sufficiently large in terms of the parameters related to the mollifiers $M(\chi,\kappa)$ defined in section \ref{sec:mollifier}.
\end{theorem}

The proof of this theorem is given in section \ref{sec:proofHighMoments}. As in \cite{DFL-2}, we start with an upper bound for $\log |L(1/2,\chi)|$ given in Lemma \ref{lem:logderLbound}. The main contribution to \eqref{momentupperbound} comes from prime ideals and their squares that appear in this upper bound. Most of section \ref{sec:HighMoments} deals with the former, which also leads to the definition of the mollifiers given in section \ref{sec:mollifier}. The contribution of squares of primes is dealt with separately in section \ref{sec:squares}.

Our next result gives an asymptotic formula for the first mollified moment of the family and is unconditional.
\begin{thm}\label{thm:firstmoment}
Let $\sF$ be the family defined in \eqref{family}. Then, 
\eqn{\label{eq:firstmoment}
\sum_{\chi \in \sF(X)} L(1/2, \chi) M(\chi,1)  = C_X X \log X  + o_X(X\log X),}
where $C_X \in (c_0,c_1)$ and  $C_X, c_0, c_1$ are given by \eqref{CX}, \eqref{c0} and \eqref{c1}, respectively.
\end{thm}
Using the definition of the mollifier $M(\chi, 1)$ in \eqref{eq:firstmoment} leads to the twisted sum \eqref{1stMoment} of $L(1/2, \chi)$ over the family of characters $\chi$. Then, by employing the approximate functional equation in \eqref{AppFncEqnatOnehalf} for each $L$-function $L(s,\chi)$ gives a character sum over the family and a dual sum of shifted cubic Gauss sums given by \eqref{S123}. The former is split into two sums $S_1$ that runs over the cubes, and the remaining terms, which constitute the sum $S_2$. The first sum, $S_1$, contributes to the main term in \eqref{eq:firstmoment} and is discussed in section \ref{sec:S1}. The second sum, $S_2$, along with the dual sum, $S_3$, gives rise to the error terms.  These error terms are estimated using Heath-Brown's cubic large sieve inequality \cite[Theorem 2]{HB2000}, a Polya-Vinogradov type inequality for $S_2$ and a bound on shifted cubic sums for $S_3$ developed in \cite{DG}. 

\section{Preliminary}

\subsection{Notation}
For $n \ge 2$, we put $\log_n X = \log (\log_{n-1} X)$, where $\log_1 X = \log X$. We use $\N\fa$ for the norm of the ideal $\fa$, and $\N(a)$ for the norm of the ideal $a\ring$. We use 
$e(z) = \exp ( 2\pi i z)$ and put $\tr(z) = z + \bar z$. We use $\fp, \fq$ for prime ideals, and write  $\fa \in \cI$ to mean that $\forall\fp, \fp \mid \fa \Rightarrow \N\fp \in \cI$ for an interval $\cI$. 

The following lemma is used to find optimal bounds in the proofs of several results throughout the paper. 
\begin{lem}[{\cite[Lemma 2.4]{GraKol}}] \label{balancing} 
Suppose that $$L(H) = \sum_{i=1}^m A_i H^{a_i} + \sum_{j=1}^n B_j H^{-b_j},$$ where
$A_i, B_j, a_i, b_j$ are positive, and that $H_1 \le H_2$. Then, there is some $H$ with $H_1 \le H \le H_2$ such that
$$
L(H) \ll \sum_{i=1}^m \sum_{j=1}^n \left(A_i^{b_j} B_j^{a_i}\right)^{1/(a_i+b_j)} + \sum_{i=1}^m A_i H_1^{a_i} + \sum_{j=1}^n B_j H_2^{-b_j},
$$
where the implied constant depends only on $m$ and $n$.
\end{lem}

\subsection{Cubic Characters}
The ring of integers $\ring$ of $K$ has class number one and six units $\left\{ \pm 1, \pm \omega, \pm \omega^2 \right\}$. Each non-trivial principal ideal $\fn$ co-prime to $3$ has a unique generator $n \equiv 1 \mod 3$.

The cubic Dirichlet characters on $\Z[\omega]$ are given by the cubic residue symbols. For each prime $\pi \in \ring$ with $\pi \equiv 1 \mod 3$, there are two primitive characters of conductor $(\pi)$; the cubic residue character $\chi_\pi$ satisfying 
\eqs{\chi_\pi(\alpha) =  \( \frac{\alpha}{\pi}\)_3 \equiv \alpha^{(\N(\pi)-1)/3} \mod \pi\ring,}
and its conjugate $\overline{\chi}_\pi= \chi_\pi^2$. 

In general, for  $n\in\ring$ with $n \equiv 1 \mod 3$, the cubic residue symbol $\chi_n$ is defined multiplicatively using the characters of prime conductor by
\eqs{\chi_n (\alpha) = \left( \frac{\alpha}{n} \right)_3 = \prod_{\pi^{v_\pi} \| n} \chi_\pi (\alpha)^{v_\pi} .}

Such a character $\chi_n$ is primitive when it is a product of characters of distinct prime conductors, i.e. either $\chi_\pi$ or $\overline \chi_\pi = \chi_{\pi}^2 = \chi_{\pi^2}$.
Moreover, $\chi_n$ is a (cubic) \emph{Hecke} character of conductor $n \ring$ if $\chi_n(\omega) = 1$. 
Since
\eqs{\( \frac{\omega}{n}\)_3 = \prod_{\pi \mid n} \omega^{v_\pi(n) (\N(\pi)-1)/3} = \omega^{\sum_{\pi \mid n} v_\pi(n) (\N(\pi)-1)/3} = \omega^{(\N(n) - 1)/3},}
we conclude that a given Dirichlet character $\chi$ is a primitive cubic Hecke character provided that $\chi = \chi_n$, where 
\begin{enumerate}
\item $n=n_1 n_2^2$, where $n_1, n_2$ are square-free and co-prime, and
\item  $\N(n) \equiv 1 \mod 9$, or equivalently, $\N(n_1) \equiv \N(n_2) \mod 9$.
\end{enumerate}
In this case, $\chi$ has conductor $n_1 n_2\ring$. 

We recall the cubic reciprocity theorem for cubic characters.
\begin{lem}[{\cite[page 114, Theorem 1]{IR}}] Let $m,n \in \ring, m,n \equiv \pm 1 \mod 3.$ Then,
$$
\( \frac{m}{n}\)_3 = \( \frac{n}{m}\)_3.
$$
\end{lem}

\subsection{The family $\sF$ of cubic Dirichlet characters}
We shall consider the family $\sF$ of cubic characters $\chi_{c_2}\overline{\chi_{c_1}}$ given by cubic residue symbols, where $c_1, c_2 \equiv 1 \mod 3\in \Z[\omega]$ are square-free and co-prime, and $c_2c_1^2 \equiv 1 \mod 9$. We naturally exclude the case $c_1=c_2=1$. By a slight abuse of notation (dropping the letter $\chi$), we shall write
\bsc{ \label{family}
\sF &= \Big\{ c_2c_1^2 \in \ring\setminus\{1\} :  \begin{array}{l} 	c_1, c_2 \equiv 1 \mod 3 \text{ both square-free},\\
(c_1,c_2)=1, \;c_2c_1^2 \equiv 1 \mod 9
\end{array} \Big\} \\
&= \Big\{ qd \in \ring\setminus\{1\} :  \begin{array}{l}
q, d \equiv 1 \mod 3, q \text{ square-free},\\
d \mid q, \;qd \equiv 1 \mod 9
\end{array} \Big\}. 
}
Note that for $c=qd \in \sF$, the conductor of the character $\chi_{qd}$ is $q\ring$. We shall write $\cond{c}$ in short, instead of $\cond{\chi_c}$, and $c \in \sF(X)$ to mean that $c\in \sF$ and $\cond c \le X$.

\subsection{Cubic Gauss Sums}
For any $n\equiv 1 \mod 3$, the shifted cubic Gauss sum is defined by 
\eqs{
g(r,n) = \sum_{\alpha \mod n} \chi_n (\alpha) \e \bigl(\tr( r \alpha/n)\bigr).}
The following two lemmas are classical results about cubic Gauss sums which can be found in \cite{HBP1979}, or easily checked, and we include them without proof.

\begin{lem} \label{cubic-GS-lemma1} 
Let $n, n_1, n_2 \equiv 1 \mod 3$ and $s, r$ be elements of $\Z[\omega]$. \\ 
If $(s, n)=1$, 
\eqs{ 
g(rs, n) = \overline{\chi_{n}}(s) g(r,n).} 
If $(n_1, n_2)=1$, 
\eqs{ 
g(r, n_1 n_2) = \overline{\chi_{n_2}(n_1)} g(r,n_1) g(r, n_2) = g(rn_1, n_2) g(r, n_1).}
\end{lem}
\begin{lem} \label{cubic-GS-lemma2}
Let $\pi \equiv 1 \mod 3$ be a prime,  $(\pi, r)=1$, where $r\equiv 1 \mod 3$. Let $k, j$ be integers with $k > 0$ and $j \ge 0$.

If $k=j+1$, 
\eqs{g(r\pi^j,\pi^k) = \N(\pi^j) \times \ar{-1 & 3\mid k\\[2mm] g(r,\pi) & k\equiv 1 \mod 3 \\[2mm] \overline{g(r,\pi)} & k \equiv 2 \mod 3.  }}

If $k \neq j+1$,
\eqs{g(r\pi^j,\pi^k) = \begin{cases}
\varphi_K(\pi^k)& \text{if } 3\mid k, \, k \le j \\ 0 & \text{otherwise.}  
\end{cases}}
\end{lem}

For any character  $\lambda$, we put $g_\lambda (r,b) = \lambda(b) g(r,b)$ and define
\eqs{
h (r,\lambda, s) = \su{b \equiv 1 \mod 3 \\ (b,r)=1}  g_\lambda(r,b) \N(b)^{-s},
}

\subsection{Hecke $L$-functions}
Let $\chi$ be a primitive Hecke character to some modulus $\fm = m\ring$ of trivial infinite type. The Hecke $L$-function associated with $\chi$ is given by
\eqs{L(s, \chi) = \sum_{\fa} \chi (\fa) (\N\fa)^{-s}}
where $\fa$ runs over the integral ideals of $\ring$ and we put
$\chi (\fa) = 0$ whenever $\fa$ is not coprime to $\fm$. The completed Hecke $L$-series is then defined by
\eqs{\Lambda (s,\chi) = (|d_K| N(m))^{s/2} (2\pi)^{-s} \Gamma(s) L(s,\chi),}
where $d_K=-3$ is the discriminant of $K$.
\begin{prop}[{\cite[VII. Cor. 8.6]{Neukirch}}]
The completed $L$-series $\Lambda (s,\chi)$ is entire, provided $\chi$ is primitive and $\fm \neq \ring$. Furthermore, it satisfies the functional equation
\eqs{\Lambda(s,\chi) = W(\chi) (\N\fm)^{-1/2} \Lambda(1-s,\overline{\chi}),}
where 
\eqs{
W(\chi)  = \su{x \mod \fm \\ (x,\fm)=1} \chi(x) \e \bigl(\tr(x/m\sqrt{-3})\bigr) }
is the associated cubic Gauss sum, and $x$ varies over a system of representatives of $\(\ring/\fm\)^\times$. 
\end{prop}

\subsection{Dedekind Zeta Function and the prime ideal theorem}
Note that 
\eqs{\zeta_K (s) 
= \zeta(s) \prod_{p \neq 3}  \Bigl( 1 - \Bigl(\frac{-3}{p}\Bigr) p^{-s} \Bigr)^{-1} = \zeta(s) L(s, (\tfrac{-3}{\cdot}) ).
}
We know by the class number formula of Dirichlet that
\eqn{\label{ZetaRes}
L(1, (\tfrac{-3}{\cdot}) ) = \frac{2\pi h(-3)}{\omega_K \sqrt{|\fd_K|}} = \frac{\pi}{3\sqrt{3}},
} 
which gives the residue of $\zeta_K(s)$ at $s=1$. 

By \cite[Theorem 2]{Serre} and \cite{Stark} it follows that there exist positive absolute constants $c_1, c_2$ such that 
\eqn{\label{PNT} 
\pi_K (x) = \li (x) + O  \bigl(x \exp ( -c_1 \sqrt{\log x} )\bigr) 
}
for $x \ge c_2$, where $\pi_K (x)$ denotes the number of prime ideals of $\ring$ whose norm does not exceed $x$.

\subsection{Family Related Results}
\begin{lem}[{\cite[Chapter 17, 1st Lemma]{Davenport}}] \label{perron}
Let
\eqs{
\delta(y) = \begin{cases}
0 & \text{if } 0 < y < 1 \\
1/2  & \text{if } y=1 \\
1 & \text{if } y > 1.
\end{cases}}
For any $y>0, a>0, T>0$,
\eqs{\bigg| \frac 1 {2\pi i} \int_{a-iT}^{a+iT} \frac{y^s}{s} ds - \delta(y) \bigg| < 
\begin{cases}
y^a \min \{ 1, T^{-1} |\log y|^{-1} \} & \text{if } y \neq 1\\
aT^{-1} & \text{if } y=1.
\end{cases}
}
\end{lem}
The next result gives the family size and is also used in establishing main term of the first moment. The error term below can be improved to $X^{2/3+\eps}$, but is not required for our purposes.
\begin{lem} \label{familysize}
For any $n \in \ring$, 
\eqs{\su{c \in \sF(X)\\ (c, n)=1}  1  = 
C_1 (n) X \log X  + C_2 (n) X + O_\eps  \Bigl( X^{3/4+\eps} \N(n)^\eps \Bigr), }
where 
\bsc{\label{C12}
C_1 (n) &= \frac {4\pi^2 F_{\psi_0}(1;n)} {2187} \\
C_2 (n) & = \frac {4F_{\psi_0}(1;n)} {81}  \Big \{\frac{\pi^2\log (3/e)}{27}   +  \lim_{s\to 1} \bigl((s-1)^2\zeta_K^2(s) \bigr)' \Big\} + \frac{4\pi^2}{2187}  F_{\psi_0}'(1;n) }
and $F_\psi (s;n)$ is defined in \eqref{Fpsi(s)}.
\end{lem}

\begin{proof}
Using Lemma \ref{perron} and proceeding similar to \cite[Chapter 17]{Davenport}, we can show for $T \ge 1$ and $\rho  = 1 + 1/(2\log X)$ that
\eqs{\su{c \in \sF(X)\\ (c, n)=1} 1 = \frac 1 {2\pi i } 
\int_{\varrho-iT}^{\varrho+iT} X^s  \biggl(\su{c \in \sF \cup \{1\}\\ (c, n)=1} \frac 1 {\cond{c}^s} \biggr) \frac{ds} s + O \Bigl( \red{X^\eps + X^{1+2\eps}T^{-1}} \Bigr).
}

Note that by the definition of the family $\sF$ (see equation \eqref{family}) we have
\bs{
\su{c \in \sF \cup \{1\}\\ (c, n)=1} \frac 1 {\cond{c}^s} &= 
\su{q \equiv 1 \mod 3\\ q \text{ square-free} \\ (q,n)=1} \frac 1 {\N(q)^s} \su{d \equiv 1 \mod 3\\ d \mid q\\ qd \equiv 1 \mod 9} 1 \\
&= \frac 1 {h_{(9)} } \sum_{\psi \mod 9} G_{\psi}(s; n),
}
where $\psi$ runs over the ray class characters modulo 9, $h_{(9)} = 9$, and for any ray class character $\Psi$ of conductor not co-prime to $3$, and for any $\alpha \in \ring$, 
\bsc{\label{G_psi}
G_\Psi(s; \alpha) 
&= \pr{\pi \equiv 1 \mod 3\\ \pi \nmid \alpha} \Bigl( 1 + \frac{\Psi(\pi)}{\N(\pi)^s} \bigl(1+\Psi(\pi)\bigr)\Bigr) = L(s,\Psi) L(s,\Psi^2) F_\Psi(s; \alpha),}
where
\eqs{L(s,\Psi) = \pr{\pi \equiv 1 \mod 3} \Bigl( 1 - \frac{\Psi(\pi)}{\N(\pi)^s} \Bigr)^{-1} = \sum_\fa \frac{\Psi(\fa)}{(\N\fa)^s} \qquad (\text{since }\Psi(1-\omega) =0)}
is the Hecke L-function associated with the ray class character $\Psi$, and
\bsc{\label{Fpsi(s)}
F_\Psi (s;\alpha) &= \pr{\pi \equiv 1 \mod 3\\ \pi \nmid \alpha}  \Bigl( 1 - \frac{\Psi^2 (\pi) + \Psi^3 (\pi)+\Psi^4(\pi)}{\N(\pi)^{2s}} + \frac{\Psi^4(\pi)+\Psi^5(\pi)}{\N(\pi)^{3s}}\Bigr) \\
&\qquad \cdot \pr{\pi \equiv 1 \mod 3\\ \pi \mid \alpha} \Bigl( 1 - \frac{\Psi(\pi)}{\N(\pi)^s} \Bigr) \Bigl( 1 - \frac{\Psi^2(\pi)}{\N(\pi)^s} \Bigr).}
Now, the integral above can be written as 
\bs{\frac 1 {2\pi i h_{(9)}} \sum_{\psi \mod 9} 
\int_{\varrho-iT}^{\varrho+iT} X^s   G_{\psi}(s; n)  \frac{ds} s.}
Note that since 
\eqs{
|F_\Psi (s;\alpha)| \le \pr{\pi \equiv 1 \mod 3\\ \pi \mid \alpha} \Bigl( 1 + \frac 1 {\N(\pi)^{\sigma}} \Bigr)^2
\exp \biggl(\su{\pi \equiv 1 \mod 3} \Bigl( \frac 3 {\N(\pi)^{2\sigma}} + \frac 2 {\N(\pi)^{3\sigma}} \Bigr) \biggr), }
$F_\Psi(s;\alpha)$ converges absolutely for any $\alpha \in \ring$ in the region $\re s \ge 1/2+\eps$, uniformly for any character $\Psi$, whereas $L(s,\Psi)$ and $L(s,\Psi^2)$ both have analytic continuation to entire functions except for a simple pole at $s=1$ if the characters are principal. Therefore, we see that $G_{\psi} (s; \alpha)$ is analytic in the region $\re s \ge 1/2 + \eps$, except for a double pole at $s=1$ when $\psi$ is principal, and a simple pole if $\psi$ has order 2. Since $h_{(9)} = 9$, there is no character of order 2. Therefore, shifting the contour to $\re s = 1/2 +  \eps$, we pick up the residue at $s=1$ when $\psi = \psi_0$, and using the crude estimate $F_\psi (s;n) \ll (\N(n))^\eps$ and the classical convexity bound $\L(s,\chi) \ll (\cond{\chi}^{1/2}(1+|t|))^{1+\eps -\sigma}$ for $-\eps \le \sigma \le 1 + \sigma$, we deduce that 
\bs{ 
\int_{1/2+\eps}^\varrho X^s G_{\psi} (s;n) \frac{ds}s &\ll \bigl(X^{1/2+\eps} + X T^{2\eps-1} \bigr) (N(n))^\eps \\
\int_{1/2+\eps - iT}^{1/2+\eps + iT} X^s G_{\psi} (s; n) \frac{ds}s  &\ll X^{1/2+\eps} T (N(n))^\eps.
}
Combining all the estimates and applying Lemma \ref{balancing} with $T \in [1,X]$, we get the claimed error term.

Since $L(s,\psi_0) L(s,\psi_0^2) = (1-3^{-s})^2 \zeta_K^2 (s)$ has a double pole at $s=1$, the residue is given by
\bs{
&\lim_{s\to 1} \bigl(X^s s^{-1} F_{\psi_0}(s;\fm) (1-3^{-s})^2 (s-1)^2 \zeta_K^2 (s) \bigr)' \\
& =\frac 49  F_{\psi_0}(1;\fm) (\mathop{\res}_{s=1} \zeta_K(s))^2 X \log X \\
&\quad + \frac 49 X \biggl(  \log (3/e) F_{\psi_0}(1;\fm) (\mathop{\res}_{s=1} \zeta_K(s))^2 \\
&\qquad + F_{\psi_0}'(1;\fm) (\mathop{\res}_{s=1} \zeta_K(s))^2 + F_{\psi_0}(1;\fm) \lim_{s\to 1} \bigl((s-1)^2\zeta_K^2(s) \bigr)' \biggr).}
The desired result follows upon using \eqref{ZetaRes}.
\end{proof}

\begin{lem} \label{charsumoverfamily}
Let $\chi$ be a non-principal cubic ray class character of conductor $\fm=m\ring$ with $(m,3)=1$. Assuming that GRH holds, the inequality 
\bs{\su{c \in \sF(X)} \chi (c) &\ll_\eps X^{1/2+\eps} n_\chi^{\eps}}
holds.
\end{lem}
\begin{proof}
We only give a sketch of the proof as it is similar to that of Lemma \ref{familysize}. Using Perron's formula, we have as before for $T \ge 1$ that
\eqs{\su{c \in \sF(X)} \chi (c) = \frac 1 {18\pi i} \sum_{\psi \mod 9} 
\int_{\varrho-iT}^{\varrho+iT} X^s   G_{\psi\chi}(s; m)  \frac{ds} s + O \Bigl( \red{X^\eps + X^{1+2\eps}T^{-1}} \Bigr),
}
where $G_\Psi (s; m)$ was defined in \eqref{G_psi}. We see that $G_{\psi\chi} (s; \alpha)$ is analytic in the region $\re s \ge 1/2 + \eps$ for each $\psi$, since $\chi$ is not principal and the conductors of $\psi$ and $\chi$ are coprime. We shift the contour to $\re s = 1/2 + \eps$. Using the estimate $F_\Psi (s;\fm) \ll (\N\fm)^\eps$ together with the bound $L(s,\chi) \ll  (|t|\cond\chi )^\eps$, we deduce that 
\bs{ \int_{1/2+\eps}^\varrho X^s G_{\psi\chi} (s;\fm) \frac{ds}s &\ll X^{1/2+\eps} T^{- 1} (n_{\psi\chi})^{3\eps} + X T^{\eps-1} (n_{\psi\chi})^{3\eps}\\
\int_{1/2+\eps - iT}^{1/2+\eps + iT} X^s G_{\psi\chi} (s;\fm) \frac{ds}s  &\ll X^{1/2+\eps} T^{\eps} (n_{\psi\chi})^{3\eps}.
}
The result follows combining all the estimates and choosing $T= X^{1/2}$. 
\end{proof}

\section{Definition of the Mollifiers} 
To define the mollifiers we shall need the following result first, the proof of which will be skipped as it is very similar to \cite[Theorem 2.1]{chandee}. 
\begin{lem} \label{lem:logderLbound}
Let $\lambda_0 \approx 0.491225$ be the unique number satisfying $e^{-\lambda_0} = \lambda_0 + \lambda_0^2/2$. Assuming GRH for $L(s,\chi)$, the inequality  
\bs{\log |L(1/2+it,\chi)| &\le \re \sum_{1 < \N\fa \le x} \frac{\Lambda(\fa) \chi (\fa)}{\N(\fa)^{1/2+\lambda/\log x+it}\log \N\fa}  \frac{\log (x/\N\fa)}{\log x} \\
&\quad + \frac{1+\lambda}{2\log x} \log \Bigl(\frac{3\cond{\chi}|1/2+\lambda/\log x + it|^2}{4\pi^2}\Bigr) \\
&\quad + \frac{2e^{-\lambda}}{x^{1/2}\log^2 x |1/2+\lambda/\log x + it|^2}}
holds for any $\lambda \ge \lambda_0$, any $t$ which does not coincide with the ordinates of zeros of $L(s,\chi)$, and any $x \ge 3$.
\end{lem}

\subsection{The Mollifier} \label{sec:mollifier}
Let $k, \kappa > 0$ be real numbers such that $k\kappa$ is an integer bounded by an absolute constant. Put $k_0 = \max\{3, 4k^2\}$. 
\begin{defn} \label{def:intervals}
Let $\alpha > 1$ and $\Theta, \beta<1$ be fixed real numbers satisfying
\eqn{\label{alphabetacond}
\red{
\alpha\beta > \frac{2k\kappa+3}{2}.}
}
Assume that $X$ satisfies 
\eqn{\label{Xcond}
(\log_2 X)^{\alpha\beta} > 2, \quad \alpha \log_3 X \ge 1 - \log \Theta.
}
Let $J = \fl{\log (\Theta (\log_2 X)^\alpha)}$, and for each $0 \le j \le J$, put 
\eqs{
\cI_j = (X^{\theta_{j-1}},  X^{\theta_j}],  \qquad \ell_j = 2 \fl{\theta_j^{-\beta}},}
where 
\eqs{
\theta_{-1} = \frac{\log k_0}{\log X}, \qquad \theta_j = \frac{e^j}{(\log_2 X)^\alpha} \quad (j \ge 0).
}
\end{defn}


Note that $\ell_j$ is a positive even integer since each $\theta_j < 1$.
\begin{lem}[{\cite[Section 3, Lemma 1]{RadSo}}] \label{exp(x)bound}
For any integer $\ell \ge 0$, put
\eqs{
E_\ell (x) =  \sum_{0 \le n \le \ell} \frac{x^n}{n!}.
}
Then, $E_{2\ell}(x) > 0$ and is concave up. Furthermore,
\bs{e^x &< E_{2\ell} (x) \qquad (x \le 0)\\
e^x &< \bigl(1+ e^{-2\ell} \bigr)E_{2\ell}(x) \qquad (x \le 2\ell/e^2).}  
\end{lem}
\begin{remark}\label{Elbound}
Note that for $x \ge \ell$, 
\eqs{E_\ell (x) \le (\ell+1) \frac{x^\ell}{\ell!}}
as $x^n/n!$ is increasing for $0 \le n \le \ell$. In particular, $E_\ell (x) < x^\ell$ if $x \ge \ell \ge 3$.
\end{remark}

\begin{defn} 
For $0 \le j \le J$, let   
\eqn{\label{Djk}
D_{j,k} (c) =  \prod_{r=0}^j (1+e^{-\ell_r}) E_{\ell_r} \bigl(k\re F_r (c, j) \bigr),}
and
\eqs{
S_{j,k} (c)  = \exp \biggl(k\re \sum_{k_0 < \N\fp \le X^{\theta_j/2}} \frac{\chi_c (\fp^2)}{2\N(\fp)^{1+2/(\theta_j\log X)}}  \Bigl(1 - \frac{2\log \N\fp}{\theta_j\log X}\Bigr) \biggr),}
where 
\bsc{\label{Fandf}	
F_r (c, j) &= \sum_{\fp \in \cI_r} \frac{\chi_c (\fp) f(\fp,j) }{\N(\fp)^{1/2}} , \quad f(\fp,j) = \frac 1 {(\N\fp)^{1/(\theta_j\log X)}} \biggl(1 - \frac{\log \N\fp}{\theta_j\log X} \biggr).}
\end{defn}
\begin{defn} 
For each $r=0, \ldots, J$, let
\eqn{\label{Tr}
\cT_r = \Bigl\{c\in \sF(X) : \max_{r\le j \le J} \re F_r (c, j) \le \frac{\ell_r}{ke^2} \Bigr\}.
}

\end{defn}

\begin{lem} \label{InT0orNot}
For each $c\in\sF(X)$, either 
$c \not\in \cT_0$ or
\bs{
|L(1/2,\chi_c)|^k &\le \su{0 \le j < J\\ j < u \le J} \exp \bigl( k/\theta_j + O(kk_0) \bigr) D_{j,k} (c) S_{j,k} (c) \Bigl( \frac{e^2 k \re F_{j+1} (c, u)}{\ell_{j+1}}\Bigr)^{s_{j+1}} \\
& \quad + \exp \bigl(k/\theta_J + O(kk_0)\bigr) D_{J,k} (c) S_{J,k} (c) 
}
for any positive even integers $s_1, \ldots, s_J$.
\end{lem}
\begin{proof}
Taking $t=0$ and $\lambda=1$ in Lemma \ref{lem:logderLbound} we have for any $x \ge 3$ and $c\in\sF$ that
\bs{
\log |L(1/2,\chi_c)| &< \re \sum_{1 < \N\fa \le x} \frac{\Lambda(\fa) \chi_c (\fa)}{\N(\fa)^{1/2+ 1/\log x }\log \N\fa}  \frac{\log (x/\N\fa)}{\log x} \\
&\quad + \frac 1 {\log x} \log \Bigl(\frac{3 \cond c (1/2+1 /\log 3)^2}{4\pi^2}\Bigr) + \frac{8}{e x^{1/2}\log^2 x}.}
The contribution of powers of prime ideals $\fp^n$ with $n \ge 3$ is
\eqs{\le \sum_{n \ge 3, (\N\fp)^n \le x} \frac 1{n(\N\fp)^{n(1/2+1/\log x)}}  \frac{\log (x/(\N\fp)^n)}{\log x} < \sum_\fp (\N\fp)^{-3/2}.}
Moreover, removing the terms, if any, with $3 \le \N\fp \le k_0$ we end up with the inequality
\bsc{\label{logLsimplified}
\log &|L(1/2,\chi_c)| 
< \re \biggl(\sum_{k_0 < \N\fp \le x} \frac{\chi_c (\fp)}{\N(\fp)^{1/2+1/\log x}}  \frac{\log (x/\N\fp)}{\log x} \biggr)\\
&\qquad + \re \biggl(\sum_{k_0 < \N\fp \le \sqrt x} \frac{\chi_c (\fp^2)}{2\N(\fp)^{1+2/\log x}} \frac{\log (x/\N\fp^2)}{\log x}\biggr) + \frac {\log \cond c} {\log x}  + O(k_0).}

Now, assume that $c\in \cT_0$. If $c \in \cT_r$ for every $0 \le r \le J$, then we choose $x = X^{\theta_J}$ in \eqref{logLsimplified}, multiply both sides by $k$ and exponentiate, then apply Lemma \ref{exp(x)bound}. This gives the last term in the claimed upper bound above. Otherwise, choose $x= X^{\theta_j}$ in \eqref{logLsimplified}, where $j+1\ge 1$ is the smallest $r$ for which $c \not\in\cT_r$. Note that $k\re F_r (c, j) \le \ell_r/e^2$ for each $0 \le r \le j$ since $c \in \cT_r$ for each such $r$, and we can apply Lemma \ref{exp(x)bound}. Furthermore, since $c \not\in \cT_{j+1}$, we have $(e^2k/\ell_{j+1})\re F_{j+1} (c, u) > 1$ for some $u=u(j)>j$. This gives the first term in the claim. 
\end{proof}

\begin{defn}[The Mollifier]
For each $c\in \sF(X)$, we define the mollifier by 
\bsc{\label{mollifier}
M(c, \kappa) &= \prod_{j=0}^J M_j (c, \kappa) \\
M_j (c, \kappa) &= E_{\ell_j} \bigl(-\tfrac 1 \kappa  F_j (c, J)\bigr).}
\end{defn}

Using the definition of $E_\ell (x)$, we see that
\bs{M_j (c, \kappa) &= \sum_{n=0}^{\ell_j} \frac 1 {n!} \biggl(-\frac 1 \kappa  \sum_{\fp \in \cI_j} \frac{\chi_c (\fp)}{\N(\fp)^{1/2}} f(\fp, J) \biggr)^n \\
&= \sum_{n=0}^{\ell_j} \frac 1 {n!} \biggl(-\frac 1 \kappa \biggr)^n \sum_{\fp_1, \ldots, \fp_n \in \cI_j} \frac{\chi_c (\fp_1\cdots \fp_n)}{\N(\fp_1\cdots \fp_n)^{1/2}} f(\fp_1\cdots \fp_n, J)  \\
&= \sum_{n=0}^{\ell_j} \biggl(-\frac 1 \kappa \biggr)^n 
\su{\fa\in \cI_j\\ \Omega(\fa) = n} \frac{\chi_c (\fa)}{\N(\fa)^{1/2}} f(\fa, J) \nu(\fa) \\
&= \su{\fa \in \cI_j\\ 0 \le \Omega(\fa) \le \ell_j} \frac{\chi_c (\fa)\lambda(\fa)f(\fa, J) \nu(\fa) }{\N(\fa)^{1/2}\kappa^{\Omega(\fa)}},
}
where we extended the definition of $f$ to be a completely multiplicative function using \eqref{Fandf}, and 
\eqs{\lambda(\fa) = (-1)^{\Omega(\fa)},\qquad \nu(\fp_1^{r_1} \cdots \fp_s^{r_s}) = \frac1 {r_1! \cdots r_s!}.}
Note that $\nu$ is multiplicative and $\lambda$ is completely multiplicative. Furthermore, for any integer $n\ge 0$,
\bs{M_j (c, \kappa)^n 
&=  \su{\fa \in \cI_j} 
\frac{\chi_c (\fa)\lambda(\fa)f(\fa, J)  }{\N(\fa)^{1/2}\kappa^{\Omega(\fa)}} \nu_n (\fa; \ell_j)
}
where
\eqs{\nu_n (\fa; \ell_j) = \su{\fa_1, \ldots, \fa_n \\
\fa_1 \cdots \fa_n = \fa\\ 0 \le \Omega(\fa_i) \le \ell_j } \nu(\fa_1) \cdots \nu(\fa_n).}
\begin{rem} \label{rem:nu}
The inequality $\nu_n (\fa; \ell_j) \le \nu_n (\fa)$ holds with equality when $\Omega(\fa) \le \ell_j$, where $\nu_n (\cdot)$ is the $n$-fold convolution of $\nu$. Furthermore, 
\eqs{\nu_n (\fp^m) = \frac{n^m}{m!}.}
\end{rem}

\section{Upper Bound For Mollified Moments} \label{sec:HighMoments}
We start with a simple preliminary result. 
\begin{lem} \label{primesumestimates}
There exists a number $X_0(\Theta) > 0$ such that for $X \ge X_0$, 
\begin{enumerate}
	\item $\su{\N\fp > k_0} \frac 1 {(\N\fp)^{3/2}} < 1$,

\item For $\sigma > 1$, $\sum_{\fp \in \cI_r} \frac 1 {(\N\fp)^\sigma} < \frac 2 {(\sigma-1)\fl{X^{\theta_{r-1}}}^{\sigma-1}}$,

\item  $\sum_{\fp \in \cI_{0}}  \frac {1} {\N\fp} < 2\log \log X^{\theta_0}$,

\item $\su{X^{\theta_j} < \N\fp \le X^{\theta_J}} \frac{1} {\N\fp} = J-j + O\(\frac 1 {\log X^{\theta_j}} \), \qquad 0 \le j \le J-1$, 

\item $\su{k_0 < \N\fp \le X^{\theta_j}} \frac{\log^2 \N\fp} {\N\fp} < 2 (\log X^{\theta_j})^2, \qquad 0 \le j \le J$.
\end{enumerate}

\end{lem}
\begin{proof}

\begin{enumerate}
\item \bs{
\su{\N\fp > k_0} \frac 1 {\N(\fp)^{3/2}} &< \su{3 < \N\fp \le 31} \frac 1 {\N(\fp)^{3/2}} + \su{\N\fp > 31 } \frac 1 {\N(\fp)^{3/2}} \\
&< 	\frac{1}{2^3}+\frac{1}{5^3}+\frac{2}{7^{3/2}}+\frac{1}{11^3}+\frac{2}{13^{3/2}}+\frac{1}{17^3} \\
&\quad +\frac{2}{19^{3/2}}+\frac{1}{23^3}+\frac{1}{29^3}+\frac{2}{31^{3/2}}+ \frac 23 < 1}
since 
\bs{
\su{\N\fp > 31 } \frac 1 {\N(\fp)^{3/2}} &< \sum_{p > 36} \frac 2 {p^{3/2}} < \sum_{n > 36} \frac 2 {n^{3/2}} < 2\int_{36}^\infty x^{-3/2} dx = \frac 2 3.
}
\item We have
\eqs{\sum_{\fp \in \cI_r} \frac 1 {(\N\fp)^\sigma} 
< 2\int_{\fl{X^{\theta_{r-1}}}}^\infty x^{-\sigma} dx = \frac 2 {(\sigma-1)\fl{X^{\theta_{r-1}}}^{\sigma-1}}.}	
\item Let $z= \sqrt{\log X^{\theta_0}} > k_0$. Then, for $X$ large enough, using \eqref{PNT} implies
\bs{
\sum_{\fp \in \cI_{0}}  \frac {1} {\N\fp} &< \sum_{k_0 < p \le z} \frac 2 p +  \sum_{z < \N\fp \le X^{\theta_0}} \frac 1 {\N\fp}\\
& < 2 \int_1^z \frac {dx} x + \log \(\frac{\log X^{\theta_0}}{\log z} \) + O\(\frac 1 {\log z}\) \\ &= 2 \log z + \log \(\frac{\log X^{\theta_0}}{\log z} \) + O \Bigl( \frac 1 {\log z} \Bigr) < 2 \log_2 X^{\theta_0}.
}
\item This follows using definition \ref{def:intervals} and \eqref{PNT}. 

\item Choose $z = \exp( (2^{-1}\log^2 X^{\theta_j})^{1/3})$. Then, it follows from \eqref{PNT} that for $X$ large enough,
\bs{\su{k_0 < \N\fp \le X^{\theta_j}} \frac{\log^2 \N\fp} {\N\fp} &\le 2 \sum_{k_0 < n \le z} \frac{\log^2 n}{n} + \su{z < \N\fp \le X^{\theta_j}} \frac{\log^2 \N\fp} {\N\fp}\\
&\le 2 \log^2 z \int_1^z \frac {dx} x + \log^2 X^{\theta_j} - \log^2 z + O \( \frac 1 {\log^3 z} \)  \\
&< 2 \log^2 X^{\theta_j}.
}
\end{enumerate}
\end{proof}

\subsection{Proof of Theorem \ref{thm:upperbound}} \label{sec:proofHighMoments}
We split the sum in \eqref{momentupperbound} into the sums 
\eqn{\label{HMsplit1}
\sum_{c \in \sF(X) \cap \cT_0} |L(1/2,\chi_c)|^k |M(c, \kappa)|^{k\kappa} }
and
\eqn{\label{HMsplit2}
\sum_{c \in \sF(X) \setminus \cT_0} |L(1/2,\chi_c)|^k |M(c, \kappa)|^{k\kappa} .
}
For each $c \in \sF(X) \setminus \cT_0$ (see \eqref{Tr}), there exists $0 \le u=u(c) \le J$ such that $\re F_0 (c, u) > \ell_0/(ke^2)$. Therefore, applying Cauchy-Schwarz inequality shows that the inequality
\eqs{
\eqref{HMsplit2}^2 \le (J+1)  \su{c \in \sF(X)} |L(1/2,\chi_c)|^{2k}  \sum_{u=0}^J  \su{c \in \sF(X)} |M(c, \kappa)|^{2k\kappa} \Bigl(\frac{ke^2\re F_0 (c, u)}{\ell_0}\Bigr)^{2s_0}
}
holds for any positive \emph{even} integer $s_0$. We choose $s_0 = (2k\kappa + 3)\ell_0$. Then, 
applying Lemma \ref{lem:matkap} with $j=-1$ together with Lemma \ref{primesumestimates} shows that the inequality 
\bs{
&\sum_{c \in \sF(X)} \bigl(\re F_0 (c, u)\bigr)^{2s_0} |M(c, \kappa)|^{2k\kappa}  \\
&\quad <   \frac {|\sF(X)| (2s_{0})! (\log_2 X^{\theta_0} )^{s_0}} {2^{s_0} (2\ell_0) !} \(\frac 5 3 \)^{2\ell_0} 
\exp \Bigl(  0.006k^3 + 2k+ 2k^2  \\
& \qquad +  \bigl(2e^{1/4}k^2 + 4ke^{1/(4k)}\bigr) \log_2 X^{\theta_0}  + \frac{0.024k^3 + k^2  + k^4}{4\fl{k_0}^2}   \\
& \qquad + 4(e-2)k^2 \( J + o_X(1)  \) \Bigr)+ o \bigl( 4^{s_0} X^{(1/2+\eps)\eta_0}\bigr)
}
holds, where $\eta_0$ is defined in \eqref{eta}. The last $o$-term can be omitted provided that
\eqs{
\eta_0 = 1 + (4k\kappa+6) \ell_0 \theta_0 + 2k\kappa \sum_{0\le r\le J} \theta_r \ell_r  \le \frac 1 {1/2+\eps} \Bigl( 1 - \frac{s_0 \log 4}{\log X}\Bigr).
}
Recalling that $\theta_r = e^r (\log_2 X)^{-\alpha}$, we see that
\eqs{\eta_0  
<  1 + 2(4k\kappa+6) \theta_0^{1-\beta} + 4k\kappa  \frac{e^{(1-\beta)} \theta_J^{1-\beta} }{e^{1-\beta}-1} .}
Since $\theta_0 < \theta_J \le \Theta$, it is sufficient to assume that 
\eqn{\label{eta0bound}
\red{
\Theta \le \Bigl(\frac{1-2\eps-o(1)}{1+2\eps} \Bigl(8k\kappa+ 12 + \frac{4k\kappa e^{1-\beta} }{e^{1-\beta}-1}\Bigr)^{-1}\Bigr)^{1/(1-\beta)}.}
}

Using the inequalities (cf. \cite[Eq. (2)]{Rob})
\eqn{\label{Stirlingbounds}
\( \frac n e\)^n \sqrt{2\pi n} \exp( 1 / (12n+1)) < n! < \( \frac n e\)^n \sqrt{2\pi n} \exp( 1 / (12n))
}
it follows that 
\eqs{\frac {(2s_{0})!} {(2\ell_0)!} \(\frac 5 3 \)^{2\ell_0} \biggl(\frac{ke^2}{\sqrt 2 \ell_0}\biggr)^{2s_0} <  \( ke(4k\kappa + 6) \(\frac {5e} {6\ell_0} \)^{2/(4k\kappa + 6)} \)^{(4k\kappa + 6)\ell_0}.
}
Since $\ell_0 > \theta_0^{-\beta}$ for $(\log_2 X)^{\alpha\beta} > 2$, which we assumed in \eqref{Xcond}, we deduce that
\bs{
&\sum_{u=0}^J \su{c \in \sF(X)} |M(c, \kappa)|^{2k\kappa} \biggl(\frac{ke^2\re F_0 (c, u)}{\ell_0}\biggr)^{2s_0} \\
&\ll J |\sF(X)| 
\exp \Biggl( \bigl( 2e^{1/4}k^2 + 4ke^{1/(4k)}\bigr) \log_2 X + 4(e-2)k^2  \( J + o_X(1)  \) + O(k) 
\\
& \quad + (4k\kappa + 6)\ell_0 \log \biggl( ke(4k\kappa + 6) \(\frac {5e} {6} \)^{1/(2k\kappa + 3)} (\log_2 X)^{1/2-\alpha\beta/(2k\kappa+3)} \biggr) \Biggr)
.
}
Combining this estimate with the result of Proposition \ref{prop:2kthmoment} and recalling that $J < \alpha \log_3 X$, it follows that
\bs{
&\eqref{HMsplit2} \ll  |\sF(X)| \exp \biggl(  k (\log_2 X)^2 + \log (\alpha \log_3 X) \\
&\quad + (4k\kappa + 6)(\log_2 X)^{\alpha\beta} \log \biggl( ke(4k\kappa + 6) \(\frac {5e} {6} \)^{1/(2k\kappa + 3)} (\log_2 X)^{1/2-\alpha\beta/(2k\kappa+3)} \biggr) \\
&\quad	+  \bigl( e^{1/4}k^2 + 2ke^{1/(4k)}\bigr) \log_2 X + 2k^2 \alpha \log_3 X \biggr) \\
& = o( |\sF(X)|),
}
as $X \to \infty$, provided that \eqref{alphabetacond} and \eqref{eta0bound} hold.

As for the sum \eqref{HMsplit1} we use the upper bound for $|L(1/2,\chi_c)|^k$ in Lemma \ref{InT0orNot} to get
\bs{
\eqref{HMsplit1} & \le  \exp \bigl(k/\theta_J + O(kk_0)\bigr) \sum_{c \in \sF(X)}  D_{J,k} (c) S_{J,k} (c) |M(c, \kappa)|^{k\kappa} \\
&\quad  + \su{0 \le j < J\\ j < u \le J} \exp \Bigl( \frac k{\theta_j} + O(kk_0) \Bigr) \\
&\qquad \cdot \sum_{c \in \sF(X)} D_{j,k} (c) S_{j,k} (c) \Bigl( \frac{e^2 k \re F_{j+1} (c, u)}{\ell_{j+1}}\Bigr)^{s_{j+1}} |M(c, \kappa)|^{k\kappa}, 
}
which holds for any positive even integers $s_1, \ldots, s_J$. To estimate the sums 
\begin{align}
\label{DJSJ}
& \sum_{c \in \sF(X)}  D_{J,k} (c) S_{J,k} (c) |M(c, \kappa)|^{k\kappa}, \\
& \label{DjSj}
\sum_{c \in \sF(X)} D_{j,k} (c) S_{j,k} (c) \Bigl( \frac{e^2 k \re F_{j+1} (c, u)}{\ell_{j+1}}\Bigr)^{s_{j+1}} |M(c, \kappa)|^{k\kappa},
\end{align}
we use Cauchy-Schwarz inequality and deduce that
\eqs{ \eqref{DJSJ}^2 \le \sum_{c \in \sF(X)}  (D_{J,k} (c) )^2  |M(c, \kappa)|^{2k\kappa} \sum_{c \in \sF(X)}  (S_{J,k} (c))^2,}
and
\eqs{ 
\eqref{DjSj}^2 \le 
\sum_{c \in \sF(X)} (D_{j,k} (c) )^2  \Bigl( \frac{e^2 k \re F_{j+1} (c, u)}{\ell_{j+1}}\Bigr)^{2s_{j+1}} |M(c, \kappa)|^{2k\kappa}  \sum_{c \in \sF(X)} (S_{j,k} (c))^2.
}
Then, we use \eqref{matkap2} in Lemma \ref{lem:matkap} together with \eqref{Sjkbound} to conclude that
\eqs{
\eqref{DJSJ}^2 \ll \cS_k  |\sF(X)|
\biggl(|\sF(X)|\exp ( 8k^2 ) \prod_{r=0}^J (1+e^{-\ell_r})^2 ( 1 + 2^{-\ell_r/2} ) + o \( X^{(1/2+\eps)\eta_J}\) \biggr), }
where $\cS_k$ is given by \eqref{Sk} and
\eqs{
\eta_J = 1  + (2k\kappa+1) \sum_{0\le r\le J} \theta_r \ell_r \le 1 + 2(2k\kappa+1) \frac{e^{(1-\beta)} \theta_J^{1-\beta} }{e^{1-\beta}-1}.
}
As before the $o$-term is $o(|\sF(X)|)$ and can be omitted under the assumption that
\eqn{\label{etaJcond}
\red{ \Theta \le \Bigl( \frac {(e^{1-\beta}-1)(1-2\eps)}{2 e^{1-\beta}  (2k\kappa+1)(1+2\eps)} \Bigr)^{1/(1-\beta)}.}
}

Finally for \eqref{DjSj}, we choose, for each $0 \le j < J$, $s_{j+1} =  2\fl{1/(4a\theta_{j+1})}$ for some $a > 1$ to be determined. Note that $s_{j+1} \ge 2k\kappa\ell_{j+1}+2$ holds for each $j$ as long as 
\eqn{\label{thetacond1}
\red{ \Theta \le \Bigl(\frac 1 {8a (k\kappa + 1)}\Bigr)^{1/(1-\beta)}.}
}
With these choices, we can then use \eqref{matkap1} in Lemma \ref{lem:matkap} together with \eqref{Sjkbound} and Lemma \ref{primesumestimates} to get for each $0 \le j < J$ that
\bs{\eqref{DjSj}^2 
&< \cS_k |\sF(X)|^2 ( 1 + o(1)) \exp ( 8k^2 )\prod_{r=0}^j (1+e^{-\ell_r})^2 ( 1 + 2^{-\ell_r/2} ) 
\\
&\quad	\cdot \Bigl( \frac{e^2 k }{2\ell_{j+1}}\Bigr)^{2s_{j+1}} \frac {(2s_{j+1})!} {\fl{\frac{2s_{j+1} - 4k\kappa\ell_{j+1}}{3}}!} \Bigl( \frac 5 3 \Bigr)^{\fl{\frac{2s_{j+1} - 4k\kappa\ell_{j+1}}{3}}} \Bigl( 1 + O \Bigl(\frac 1 {\log X^{\theta_j}} \Bigr) \Bigr)^{s_{j+1}} \\
&\quad \cdot 
\exp \biggl( e^{1/4}k^2 + 2ke^{1/(4k)} + 4(e-2)k^2 (J-j-1) + O_k  \Bigl(\frac 1 {\log X^{\theta_j} } \Bigr)\biggr),
}
where we replaced the $o$-term in \eqref{matkap1} by $o(|\sF(X)|)$ imposing the condition that 
\eqs{
\eta_j \le (1/2+\eps)^{-1} \Bigl(1 - \frac{\log 4}{2a\theta_{j+1} \log X} \Bigr),
}
and this holds provided that 
\eqn{\label{etajcond}
\red{ \Theta  \le  \Bigl(\frac{e^{1-\beta}-1} {2(2k\kappa+1)e^{1-\beta}} \( \frac{1-2\eps-o(1)}{1+2\eps} - \frac 1 a\)\Bigr)^{1/(1-\beta)}}.}

Now, write $2s_{j+1} = A$, and $\fl{(2s_{j+1}-4k\kappa\ell_{j+1})/3} = B$. Then, $B/A < 1/3$ and 
\eqn{\label{betacond}
\beta + B/A - 1 > \beta - \frac 23 - \frac {2a(8k\kappa + 1) \Theta^{1-\beta}} 3. }
Thus, $\beta + B/A - 1 > 0$ provided that
\eqn{\label{ThetaBetacond} 
\red{\beta > \frac 2 3, \qquad \Theta <  \Bigl( \frac{3\beta-2}{2a(8k\kappa+3)}\Bigr)^{1/(1-\beta)}.}
}
Assuming \eqref{ThetaBetacond} holds, it follows using \eqref{Stirlingbounds} that
\bsc{\label{magic} 
\Bigl( \frac 5 3 \Bigr)^B & \Bigl( \frac{e^2 k }{2\ell_{j+1}}\Bigr)^A \frac{A!}{B!} \\
& < \Bigl( \Bigl( \frac 5 3 \Bigr)^{B/A} \frac{e^2 k }{2\theta_{j+1}^{-\beta}}\Bigr)^A A^{A-B} \exp\( \log 2 + \frac 1 {12A} + B \log (4e)- A \) \\
& < \Bigl( \Bigl( \frac 5 3 \Bigr)^{1/3}\frac{e^2 k \theta_{j+1}^{\beta + B/A-1}}{2a^{1-B/A}}\Bigr)^A  \\
&< \exp \( \frac 1 {2a\theta_{j+1}} \log \Bigl( \Bigl( \frac 5 {3a^2} \Bigr)^{1/3} \frac{e^2 k \theta_{j+1}^{\beta + B/A-1}}{2}\Bigr) \), }
where we used $\ell_{j+1} > \theta_{j+1}^{-\beta}$ for the first inequality, which holds assuming \eqref{ThetaBetacond}, the second inequality follows since $B/A<1/3, A \le 1/(a\theta_{j+1})$, and $A \ge 4$ implies $1/(12A) + \log 2 + B \log (4e) < A$, and the third inequality follows since $B/A < 1/3$ and $A > 1/(2a\theta_{j+1})$ assuming that
\eqn{\label{thetacond4}
\red{\Theta   < \Bigl(\frac 2 {e^2 k} \Bigl( \frac {3a^2} 5 \Bigr)^{1/3}\Bigr)^{1/(\beta -2/3)}
}.
} 
Therefore, using \eqref{magic} we see that
\bs{\eqref{DjSj}^2 
&< |\sF(X)|^2 ( 1 + o(1)) \cS_k \exp ( 8k^2 +e^{1/4}k^2 + 2ke^{1/(4k)} )\prod_{r=0}^j (1+e^{-\ell_r})^2 ( 1 + 2^{-\ell_r/2} ) 
\\
&\quad \cdot \exp \Biggl( \frac 1 {2a\theta_{j+1}} \log \Bigl( \Bigl( \frac 5 {3a^2} \Bigr)^{1/3} \frac{e^2 k \theta_{j+1}^{\beta + B/A-1}}{2}\Bigr) + O \Bigl( \frac 1 {\theta_{j+1}\log X^{\theta_j}} \Bigr) \\
&\qquad \qquad + 4(e-2)k^2 (J-j-1) + O_k  \Bigl(\frac 1 {\log X^{\theta_j} } \Bigr)\Biggr). 
}

Combining the estimates for \eqref{DJSJ} and \eqref{DjSj} we conclude that 
\bs{
\eqref{HMsplit1} & <  |\sF(X)| ( 1 + o_X(1)) \cD_k 
\( C_J \exp \bigl( e^{1/4}k^2/2 + ke^{1/(4k)} \bigr)  	
+ \exp \bigl(ek/\Theta \bigr) \)
}
where
\eqs{
\cD_k = \cS_k^{1/2}
\exp \bigl( 4k^2 + O(kk_0) \bigr) \prod_{r=0}^J (1+e^{-\ell_r}) ( 1 + 2^{-\ell_r/2} )^{1/2},}
and
\bs{
C_J &= \su{0 \le j < J\\ j < u \le J} \exp \biggl( \frac k{\theta_j}  + \frac 1 {4a\theta_{j+1}} \log \Bigl( \Bigl( \frac 5 {3a^2} \Bigr)^{1/3} \frac{e^2 k \theta_{j+1}^{\beta + B/A-1}}{2}\Bigr) + O \Bigl( \frac 1 {\theta_{j+1}\log X^{\theta_j}} \Bigr) \\
&\qquad + 2(e-2)k^2 (J-j-1) \biggr) \\
&= \su{0 \le j \le  J-1} (J-j)\exp \biggl( \frac { ke^{J-j}}{\theta_J}  + \frac { e^{J-j-1}}{4a\theta_J} \log \Bigl( \Bigl( \frac 5 {3a^2} \Bigr)^{1/3} \frac{e^2 k (e^{j+1-J}\theta_J)^{\beta + B/A-1}}{2}\Bigr) \\
&\qquad + O \Bigl( \frac 1 {e^{j+1-J} e^{j-J}\theta_J^2 \log X} \Bigr) + 2(e-2)k^2 (J-1-j) \biggr) \\
&= \su{0 \le  r \le J-1} (r+1)\exp \biggl( \frac { ke^{r+1}}{\theta_J}  + \frac { e^r}{4a\theta_J} \log \Bigl( \Bigl( \frac 5 {3a^2} \Bigr)^{1/3} \frac{e^2 k (e^{-r}\theta_J)^{\beta + B/A-1}}{2}\Bigr) \\
&\qquad + O \Bigl( \frac {e^{2r}} {\theta_J^2 \log X} \Bigr) + 2(e-2)k^2r \biggr) \\
&< \su{0 \le  r \le J-1} (r+1) \exp \biggl( f(r) + O \Bigl( \frac {(\log_2 X)^{2\alpha}} {\log X} \Bigr)  \biggr).
}
Here,
\bsc{\label{R1R2}
f(r) &= R_1 e^r - R_2 r e^r + 2k^2r,\\ 
R_1 & = \frac e \Theta \biggl( ke   + \frac {1}{4a} \log \Bigl( \Bigl( \frac 5 {3a^2} \Bigr)^{1/3} \frac{e^2 k} 2 \Bigr) \\
&\qquad\quad + \frac {\log \Theta} {4a} \Bigl(\beta - \frac 23 - \frac {2a(8k\kappa + 1) \Theta^{1-\beta}} 3\Bigr)  \biggr), \\
R_2 & = \frac e {4a\Theta} \Bigl(\beta - \frac 23 - \frac {2a(8k\kappa + 1) \Theta^{1-\beta}} 3 \Bigr).
}

To finish the proof of Theorem \ref{thm:upperbound}, we observe that it is possible to choose the parameters $\eps, a, \beta, \alpha$ in the given order (for example, one possible choice, is to take
$\eps = 1/64$, $a = 6$, $\beta = 3/4$, $\alpha = 3 + 2k\kappa$),
and then choose $\Theta$ small enough in such a way that $R_1 < 0, R_2 > 0$, and \eqref{alphabetacond}, \eqref{eta0bound}, \eqref{etaJcond}, \eqref{thetacond1}, \eqref{etajcond}, \eqref{betacond}, \eqref{ThetaBetacond}, \eqref{thetacond4}, and \eqref{thetacond5} are all satisfied. This shows that $C_J \ll 1$ and for $X$ sufficiently large in terms of the chosen parameters, \eqref{momentupperbound} holds. This finishes the proof.

\subsection{Behind the Scenes}
In this part, we provide the technical results used in the proof of Theorem \ref{thm:upperbound}.

\begin{lem} \label{lem:matkap}
Let $j, u$ be integers satisfying $-1 \le j \le J-1$ and $j+1 \le u \le J$. For any positive integer $s_{j+1}$ satisfying $s_{j+1} \ge 2k\kappa\ell_{j+1}+2$,
\bsc{
\label{matkap1}
&\sum_{c \in \sF(X)} D_{j,k} (c)^2 \bigl(\re F_{j+1} (c, u)\bigr)^{2s_{j+1}} |M(c, \kappa)|^{2k\kappa}  \\
&< \exp( 8k^2 )\prod_{r=0}^j (1+e^{-\ell_r})^2 ( 1 + 2^{-\ell_r/2} ) 
\\
&\quad	\cdot  \frac {|\sF(X)| (2s_{j+1})!} {4^{s_{j+1}}\fl{\frac{2s_{j+1} - 4k\kappa\ell_{j+1}}{3}}!} \biggl(\su{\fp \in \cI_{j+1}} \frac 1 {N\fp}  \biggr)^{s_{j+1}}  \Bigl(\frac 5 3 \Bigr)^{\fl{\frac{2s_{j+1} - 4k\kappa\ell_{j+1}}{3}}} \\
&\quad \cdot \exp \biggl( \frac{0.024k^3 + k^2  + k^4}{4\fl{X^{\theta_j}}^2}   
+ \sum_{\fp \in \cI_{j+1}}  \frac{0.006k^3 + 2k+ 2k^2 }{(\N\fp)^{3/2}}
+    \frac {e^{1/4}k^2 + 2ke^{1/(4k)}} {\N\fp}
\biggr)
\\
& \qquad + \su{X^{\theta_{j+1}} < \N\fp \le X^{\theta_J}}  
\frac{4(e-2)k^2} {\N\fp}  \biggr) 
+ o \( 4^{s_{j+1}} X^{(1/2+\eps)\eta_j}\)
}
where the first line of the upper bound is to be omitted for $j=-1$. Furthermore, 
\bsc{\label{matkap2}
&\sum_{c \in \sF(X)}  (D_{J,k} (c) )^2  |M(c, \kappa)|^{2k\kappa} \\
&\quad <
|\sF(X)| \exp ( 8k^2) \prod_{r=0}^J (1+e^{-\ell_r})^2 ( 1 + 2^{-\ell_r/2} ) + o \( X^{(1/2+\eps)\eta_J}\)
}	
where
\eqn{\label{eta}
\eta_j = 1 + 2s_{j+1}\theta_{j+1} + 2k\kappa \sum_{0\le r\le J} \theta_r \ell_r + \sum_{0\le r\le j} \theta_r \ell_r .
}
and the second term $2s_{j+1}\theta_{j+1}$ is to be omitted when $j=J$.
\end{lem}
\begin{proof}
Using the definitions \eqref{Djk} and \eqref{mollifier}, we see that 
\eqs{
\sum_{c \in \sF(X)} D_{j,k} (c)^2 \bigl(\re F_{j+1} (c, u)\bigr)^{2s_{j+1}} |M(c, \kappa)|^{2k\kappa} =  \prod_{r=0}^j (1+e^{-\ell_r})^2 \sum_{c \in \sF(X)} \prod_{r=0}^J \cE(r,c),
}
where, for any $-1 \le j \le J$, 
\bs{
\cE(r,c) = \begin{cases}
\(E_{\ell_r} \bigl(k\re F_r (c, j) \bigr)\)^2 |M_r (c, \kappa)|^{2k\kappa} & \text{if } 0 \le r \le j \\
\bigl(\re F_r (c, u)\bigr)^{2s_r} |M_r (c,\kappa)|^{2k\kappa} & \text{if } r= j+1 \\
|M_r (c, \kappa)|^{2k\kappa} & \text{if } j+1 < r \le J.
\end{cases}
}
Now writing $\re z = \frac 1 2 (z+\bar z)$, we obtain for any $n \ge 0$ that
\bsc {\label{realpartofFIr(c,j)}
(\re F_r ( c,u ) )^n 
&= \frac {n!} {2^n} \su{\fa, \fb \in \cI_r \\ \Omega (\fa\fb) = n } 
\frac {f(\fa\fb, u)\nu(\fa) \nu(\fb)}{\N(\fa\fb)^{1/2}}  \chi_c (\fa\fb^2)\\
E_{\ell_r} \bigl(k\re F_r (c, j)) 
&= \su{\fa, \fb \in \cI_r \\ 0 \le \Omega (\fa\fb) \le \ell_r } (k/2)^{\Omega(\fa\fb) } 
\frac {f(\fa\fb, j) \nu(\fa) \nu(\fb)}{\N(\fa\fb)^{1/2}}  \chi_c (\fa\fb^2).
}
Therefore, 
\bsc{\label{Djksquare}
D_{j,k} (c)^2	&= \prod_{r=0}^j (1+e^{-\ell_r})^2 \su{\fa_{1r}, \fb_{1r}, \fa_{2r}, \fb_{2r} \in \cI_r \\ 
0 \le \Omega (\fa_{1r}\fb_{1r}) \le \ell_r \\
0 \le \Omega (\fa_{2r}\fb_{2r}) \le \ell_r } 
(k/2)^{\Omega(\fa_{1r}\fb_{1r}\fa_{2r}\fb_{2r}) } \\
&\qquad \cdot
\frac {f(\fa_{1r}\fb_{1r}\fa_{2r}\fb_{2r}, j) \nu(\fa_{1r}) \nu(\fa_{2r}) \nu(\fb_{1r}) \nu(\fb_{2r})}
{\N(\fa_{1r}\fb_{1r}\fa_{2r}\fb_{2r})^{1/2}}  \chi_c (\fa_{1r}\fb_{1r}^2\fa_{2r}\fb_{2r}^2).
}
Moreover, since $2k\kappa$ is a positive even integer, 
\bsc{\label{Mjto2kkappa}
|M_r (c, \kappa)|^{2k\kappa} &= \su{\fa, \fb \in \cI_r} 
\frac{\chi_c (\fa\fb^2) \lambda(\fa\fb) f(\fa\fb, J)} {\N(\fa\fb)^{1/2}\kappa^{\Omega(\fa\fb)}} \nu_{k\kappa} (\fa; \ell_r) \nu_{k\kappa} (\fb; \ell_r).
}
Using \eqref{realpartofFIr(c,j)}, \eqref{Djksquare} and \eqref{Mjto2kkappa} we see that for $0 \le r \le j$,
\bs{
\cE(r,c) &=\mathop{\sum \ldots \ldots \sum}_{
\substack{\fa_{1r},\fb_{1r}, \fa_{2r},\fb_{2r}, \fa_{3r},\fb_{3r} \in \cI_r \\ 
0 \le \Omega (\fa_{1r}\fb_{1r}) \le \ell_r\\
0 \le \Omega (\fa_{2r}\fb_{2r}) \le \ell_r \\ (c,\fm_r)=1
}} 
\frac{ \lambda(\fa_{3r}\fb_{3r}) f(\fa_{1r}\fb_{1r} \fa_{2r} \fb_{2r}, j)  f(\fa_{3r}\fb_{3r}, J)}
{\N(\fa_{1r}\fb_{1r}\fa_{2r}\fb_{2r} \fa_{3r}\fb_{3r})^{1/2}\kappa^{\Omega(\fa_{3r}\fb_{3r})}}  
\\
&\quad \cdot 
\Bigl(\frac k 2\Bigr)^{\Omega (\fa_{1r} \fb_{1r} \fa_{2r} \fb_{2r})}  \nu(\fa_{1r})\nu(\fb_{1r}) \nu(\fa_{2r})\nu(\fb_{2r})
\nu_{k\kappa} (\fa_{3r}; \ell_r) \nu_{k\kappa} (\fb_{3r}; \ell_r) 	\chi_c (\fm_r),
}
for $r=j+1$, 
\bs {\cE(r,c) &= \frac {(2s_{r})!} {4^{s_{r}}} 
\mathop{\sum \ldots \ldots \sum}_{\substack{\fa_{1r}, \fb_{1r}, \fa_{2r}, \fb_{2r} \in \cI_{r}\\ \Omega (\fa_{1r}\fb_{1r}) = 2s_{r} \\ (c,\fm_r)=1
} }
\frac {\lambda(\fa_{2r}\fb_{2r}) f(\fa_{1r}\fb_{1r}, u) f(\fa_{2r}\fb_{2r},J) 
}
{\N(\fa_{1r}\fb_{1r} \fa_{2r}\fb_{2r} )^{1/2}
\kappa^{\Omega(\fa_{2r}\fb_{2r})} 	} \\
&\qquad \cdot  
\nu(\fa_{1r})\nu(\fb_{1r}) \nu_{k\kappa} (\fa_{2r}; \ell_{r}) \nu_{k\kappa} (\fb_{2r}; \ell_{r}) \chi_c (\fm_r),}
and for $j+1 < r \le J$,
\bs{\cE(r,c) &=\su{\fa_{1r}, \fb_{1r} \in \cI_r\\ (c,\fm_r)=1} \frac{\lambda(\fa_{1r}\fb_{1r}) f(\fa_{1r}\fb_{1r},J)}
{\N(\fa_{1r}\fb_{1r})^{1/2}\kappa^{\Omega(\fa_{1r}\fb_{1r})}}
\nu_{k\kappa} (\fa_{1r}; \ell_r) \nu_{k\kappa} (\fb_{1r}; \ell_r) \chi_c (\fm_r).}
Here,
\bs{
\fm_r = \begin{cases}
\fa_{1r}\fb_{1r}^2\fa_{2r}\fb_{2r}^2 \fa_{3r}\fb_{3r}^2 & \text{if } 0 \le r \le j\\
\fa_{1r}\fb_{1r}^2\fa_{2r}\fb_{2r}^2 & \text{if } r=j+1\\
\fa_{1r}\fb_{1r}^2 & \text{if } j+1 < r \le J.	
\end{cases}
}
Note that $\fm_r$'s are pairwise co-prime since the intervals $\cI_r$ are disjoint. For each $r=0, \ldots, J$, let $\tilde \cE(r,c)$ denote the sum of the terms of $\cE(r,c)$ for which $\fm_r = \cube$, and put
\eqs{\fm = \prod_{r=0}^J \fm_r.} 
Notice that $\fm = \cube$ if and only if $\fm_r = \cube$ for each $r$. To estimate 
\eqs{\sum_{c \in \sF(X)} \(\prod_{r=0}^J  \cE(r,c) - \prod_{r=0}^J \tilde \cE(r,c)\)}
we interchange the summation over $c$ after expanding the products and taking the difference so that the innermost sum is
\eqs{\su{c \in \sF(X) \\ (c, \fm) = 1} \chi_c (\fm).}
Then, applying Lemma \ref{charsumoverfamily} yields the bound $X^{1/2+\eps} (\N\fm)^\eps$ for this sum, which in turn gives that
\bsc{\label{trashnoncubes}
&\prod_{r=0}^j (1+e^{-\ell_r})^2 \su{c \in \sF(X)} \(\prod_{r=0}^J  \cE(r,c) - \prod_{r=0}^J \tilde \cE(r,c)	\) \\
&\ll X^{1/2+\eps} \prod_{r=0}^j (1+e^{-\ell_r})^2 \biggl( E_{\ell_r}  \Bigl(k \sum_{\fp \in \cI_r} \frac {f(\fp, j)}{\N(\fp)^{1/2-\eps}}\Bigr) \biggr)^2 \\
&\cdot \biggl( \sum_{\fp \in \cI_{j+1}} \frac {f(\fp,u)}{\N(\fp)^{1/2-\eps}}  \biggr)^{2s_{j+1}}
\Biggl( \prod_{r=0}^J E_{\ell_r}\biggl(\frac 1 \kappa  \sum_{\fp \in \cI_r} \frac{f(\fp, J)}{\N(\fp)^{1/2-\eps}} \biggr) \Biggr)^{2k\kappa}
}
where the sum over the interval $\cI_{j+1}$ is to be omitted when $j = J$. 
Since
\eqs{\sum_{\fp \in \cI_r} f(\fp, j)(\N\fp)^{\eps-1/2} < 2X^{\theta_r (1/2+\eps)}} 
for any $0\le r, j \le J$, and
\eqs{2kX^{\theta_r (1/2+\eps)} > \ell_r, \qquad 2\kappa^{-1}X^{\theta_r (1/2+\eps)} > \ell_r}
for sufficiently large $X$ depending on the choice of $\alpha$ and $\beta$, it follows from Remark \ref{Elbound} 
that 
\eqn{\label{errorfigaro}
\eqref{trashnoncubes} \ll 4^{s_{j+1}} X^{(1/2+\eps)\eta_j} \prod_{r=0}^j \Bigl(\frac{(2k)^{\ell_r}(\ell_r + 1)(1+e^{-\ell_r})}{(\ell_r)!}\Bigr)^2  \prod_{r=0}^J \Bigl(    \frac{2^{\ell_r}(\ell_r+1)}{\kappa^{\ell_r}(\ell_r)!}  \Bigr)^{2k\kappa},   }
where $\eta_j$ is given by \eqref{eta}. Using Stirling's formula and the definition of $\ell_r$, it follows that $\eqref{errorfigaro} = o_X ( 4^{s_{j+1}} X^{(1/2+\eps)\eta_j})$.

We have shown so far that
\mult{\prod_{r=0}^j (1+e^{-\ell_r})^2 \sum_{c \in \sF(X)} \prod_{r=0}^J \cE(r,c) = 
\prod_{r=0}^j (1+e^{-\ell_r})^2 \sum_{c \in \sF(X)} \prod_{r=0}^J \tilde \cE(r,c) \\
+ o \( 4^{s_{j+1}} X^{(1/2+\eps)\eta_j}\)
}
for any $-1 \le j \le J$. Combining the estimates for each $\tilde \cE(r,c)$ that will be derived in the next section in \eqref{E(r)Productbound}, \eqref{E(r,c)bound2} and \eqref{E(r,c)bound3}  gives the claimed bound in \eqref{matkap1}. Similarly, using \eqref{E(r)Productbound} yields \eqref{matkap2}.
\end{proof}

\subsection{Upper Bound for $\tilde \cE(r,c)$ for $0 \le r \le j$}
Using Remark \ref{rem:nu} we can write $\tilde \cE(r,c) = \tilde \cE_1 (r,c) + \tilde \cE_2(r,c)$, where
\bs{
&\tilde \cE_1 (r,c) =
\mathop{\sum \ldots \ldots \sum}_{
\substack{\fa_{1r},\fb_{1r}, \fa_{2r},\fb_{2r}, \fa_{3r},\fb_{3r} \in \cI_r \\ \fm_r = \cube \\ (c,\fm_r) = 1
}} 
\frac{ \lambda(\fa_{3r}\fb_{3r}) f(\fa_{1r}\fb_{1r} \fa_{2r} \fb_{2r}, j)  f(\fa_{3r}\fb_{3r}, J)}
{\N(\fa_{1r}\fb_{1r}\fa_{2r}\fb_{2r} \fa_{3r}\fb_{3r})^{1/2}\kappa^{\Omega(\fa_{3r}\fb_{3r})}}  
\\
&\qquad \cdot 
\Bigl(\frac k 2\Bigr)^{\Omega (\fa_{1r} \fb_{1r} \fa_{2r} \fb_{2r})} 
\nu(\fa_{1r})\nu(\fb_{1r}) \nu(\fa_{2r})\nu(\fb_{2r}) \nu_{k\kappa} (\fa_{3r}) \nu_{k\kappa} (\fb_{3r}) 
}
and
\bs{
&\tilde \cE_2 (r,c)=
\mathop{\sum \ldots \ldots \sum}_{\substack{
\fa_{1r},\fb_{1r}, \fa_{2r},\fb_{2r}, \fa_{3r},\fb_{3r} \in \cI_r  \\ 
0 \le \Omega (\fa_{1r}\fb_{1r}) \le \ell_r\\
0 \le \Omega (\fa_{2r}\fb_{2r}) \le \ell_r\\
\max \{\Omega (\fa_{3r}), \Omega(\fb_{3r}) \} > \ell_r \\
\fm_r = \cube, (c,\fm_r) = 1
}}  
\frac{ \lambda(\fa_{3r}\fb_{3r}) f(\fa_{1r}\fb_{1r} \fa_{2r} \fb_{2r}, j)  f(\fa_{3r}\fb_{3r}, J) }
{\N(\fa_{1r}\fb_{1r}\fa_{2r}\fb_{2r} \fa_{3r}\fb_{3r})^{1/2}\kappa^{\Omega(\fa_{3r}\fb_{3r})}}  
\\
&\qquad \cdot 
\Bigl(\frac k 2\Bigr)^{\Omega (\fa_{1r} \fb_{1r} \fa_{2r} \fb_{2r})}
\nu(\fa_{1r})\nu(\fb_{1r}) \nu(\fa_{2r})\nu(\fb_{2r})
\nu_{k\kappa} (\fa_{3r};\ell_r) \nu_{k\kappa} (\fb_{3r};\ell_r) \\
& - \mathop{\sum \ldots \ldots \sum}_{\substack{
\fa_{1r},\fb_{1r}, \fa_{2r},\fb_{2r}, \fa_{3r},\fb_{3r} \in \cI_r  \\ 
\max \{\Omega (\fa_{1r}\fb_{1r}), \Omega (\fa_{2r}\fb_{2r}),\\
\Omega (\fa_{3r}), \Omega(\fb_{3r}) \} > \ell_r\\
\fm_r  = \cube, \\ (c,\fm_r) = 1
}} 
\frac{ \lambda(\fa_{3r}\fb_{3r}) f(\fa_{1r}\fb_{1r} \fa_{2r} \fb_{2r}, j)  f(\fa_{3r}\fb_{3r}, J)}
{\N(\fa_{1r}\fb_{1r}\fa_{2r}\fb_{2r} \fa_{3r}\fb_{3r})^{1/2}\kappa^{\Omega(\fa_{3r}\fb_{3r})}}  
\\
&\qquad \cdot 
\Bigl(\frac k 2\Bigr)^{\Omega (\fa_{1r} \fb_{1r} \fa_{2r} \fb_{2r})} 
\nu(\fa_{1r})\nu(\fb_{1r}) \nu(\fa_{2r})\nu(\fb_{2r})
\nu_{k\kappa} (\fa_{3r}) \nu_{k\kappa} (\fb_{3r}).
}
\subsubsection{Estimate of $\tilde \cE_1(r,c)$} 
We can write 
\bs{
\tilde \cE_1(r,c)
&= 	\mathop{\sum \ldots \ldots \sum}_{\substack{ 			\fa_{1r},\fb_{1r}, \fa_{2r},\fb_{2r} \in \cI_r \\ 			\fa_{1r} \fa_{2r} \fb_{1r}^2 \fb_{2r}^2 = \cube 	}} 	\Bigl(\frac k 2\Bigr)^{\Omega (\fa_{1r} \fb_{1r})} 	\frac{ \lambda(\fa_{2r}\fb_{2r}) f(\fa_{1r}\fb_{1r}, j)  f(\fa_{2r}\fb_{2r}, J)} 	{\N(\fa_{1r}\fb_{1r}\fa_{2r}\fb_{2r})^{1/2} \kappa^{\Omega(\fa_{2r}\fb_{2r})}}  	\\ 	&\qquad \cdot \nu_2(\fa_{1r}) \nu_2 (\fb_{1r}) \nu_{k\kappa} (\fa_{2r}) \nu_{k\kappa} (\fb_{2r}) \\
&= \mathop{\sum  \sum}_{\substack{
\fa_r,\fb_r \in \cI_r \\ 
\fa_r \fb_r^2  = \cube\\ (c,\fa_r\fb_r) =1 	}}
\frac{W(\fa_r) W(\fb_r)}{\N(\fa_r\fb_r)^{1/2}} 
= \su{\fc_r \in \cI_r\\ (c,\fc_r) =1} \frac 1 {\N(\fc_r)^{1/2}} 
\su{\fa_r,\fb_r \in \cI_r \\ \fa_r \fb_r = \fc_r \\ 
\fc_r \fb_r  = \cube 	} W(\fa_r) W(\fb_r)
}
where $W = U \ast V$ is the convolution of $U$ and $V$, and 
\eqs{U (\fa) = (k/2)^{\Omega (\fa)} 
f(\fa, j) \nu_2(\fa), \qquad V(\fa) = 
\frac{ \lambda(\fa) f(\fa, J) \nu_{k\kappa} (\fa)}
{\kappa^{\Omega(\fa)}}.}
Since all factors are multiplicative in $\fc_r$, we obtain
\eqs{
\tilde \cE_1(r,c) = \pr{\fp \in \cI_r\\ \fp \nmid c} 
\biggl(1 + \sum_{n \ge 1} \frac 1 {\N(\fp)^{n/2}} 
\mathop{\sum  \sum}_{\substack{
\fa_r,\fb_r \in \cI_r \\ \fa_r \fb_r = \fp^n \\ 
\fp^n \fb_r  = \cube 	}} W(\fa_r) W(\fb_r)\biggr). 
}
Noting that 
\bs{
W (\fp^i) 
&= \frac{k^i}{i!} \sum_{0 \le a \le i}  
\binom{i}{a} f(\fp, j)^a(-f(\fp, J))^{i-a}
= \frac{k^i}{i!} \bigl(f(\fp, j) - f(\fp, J)\bigr)^i
}
we have that 
\eqs{
\mathop{\sum  \sum}_{\substack{
\fa_r,\fb_r \in \cI_r \\ \fa_r \fb_r = \fp^n \\ 
\fp^n \fb_r  = \cube 	}} W(\fa_r) W(\fb_r) = (-1)^n b_n (\fp),}
where 
\eqs{ b_n (\fp) = \su{0 \le m \le n \\ 3 \mid n+m } W(\fp^{n-m}) W(\fp^m)  \\
= k^n \bigl(f(\fp, J) - f(\fp, j)\bigr)^n \su{0 \le m \le n \\ 3 \mid n+m }  \frac{1 }{ m! (n-m)!}.
}
Since $b_1 (\fp) = 0$, we conclude that
\eqs{
\tilde \cE_1(r,c)  = \pr{\fp \in \cI_r\\ \fp \nmid c} 
\biggl(1 + \sum_{n \ge 2} \frac {(-1)^n b_n (\fp)} {\N(\fp)^{n/2}}\biggr).
}
We will show that the sum over $n$ is an alternating series and is thus bounded above by $b_2 (\fp) /\N\fp$. First note that if $2 \mid n$, then
\bs{
&\su{0 \le m \le n+1 \\ 3 \mid n+1+m }  \frac{1 }{ m! (n+1-m)!}
= \su{0 \le m \le n/2 - 1 \\ 3 \mid n+1+m }  \frac{2}{ m! (n+1-m)!} \\ 
&< \su{1 \le m \le n/2 \\ 3 \mid n+m }  \frac{1 }{ m! (n-m)!} \le \su{0 \le m \le n \\ 3 \mid n+m }  \frac{1 }{ m! (n-m)!}.
}
If $2 \nmid n$, then $n \ge 3$ and we have
\bs{
&\su{0 \le m \le n+1 \\ 3 \mid n+1+m }  \frac{1 }{ m! (n+1-m)!}
= \su{0 \le m \le (n-5)/2 \\ 3 \mid n+1+m }  \frac{2}{ m! (n+1-m)!} + \frac{1}{ \(\frac{n+1}2\)! \(\frac{n+1}2\)!}\\ 
&< \su{0 \le m \le (n-3)/2 \\ 3 \mid n+m }  \frac{2 }{ m! (n-m)!} \le \su{0 \le m \le n \\ 3 \mid n+m }  \frac{1 }{ m! (n-m)!}.
}
Furthermore, since $e^{-x} > 1 - x$ for $x \neq 0$, we have
\bs{0 < &f(\fp,J) - f(\fp,j) = f(\fp,J) - \frac 1 {(\N\fp)^{1/(\theta_j\log X)}} \biggl(1 - \frac{\log \N\fp}{\theta_j\log X} \biggr) \\
&< 1 - \frac 1 {(\N\fp)^{1/(\theta_j\log X)}} \biggl(1 - \frac{\log \N\fp}{\theta_j\log X} \biggr) < 1 - \biggl(1 - \frac{\log \N\fp}{\theta_j\log X} \biggr)^2\\
&< \min \bigg\{ 1 , \frac{2\log \N\fp}{\theta_j \log X} \bigg \}. 
}
Hence, $0 < b_n (\fp) < (2k)^n/n! \to 0$ as $n \to \infty$ and $b_n(\fp)/(\N\fp)^{n/2}$ is decreasing since 
\eqs{\frac{b_{n+1}(\fp) (\N\fp)^{n/2} }{b_n(\fp) (\N\fp)^{(n+1)/2}} < \frac{k\bigl(f(\fp, J) - f(\fp, j)\bigr)}{ \N(\fp)^{1/2} }   
< \frac{k}{\sqrt{\N\fp}} < 1
}
for all $\fp \in \cI_r$. This shows that we have an alternating series. Therefore, 
\bsc{\label{E1(r,c)bound}
\tilde \cE_1(r,c) &<  \pr{\fp \in \cI_r\\ \fp \nmid c} \Bigl(1 + \frac{b_2(\fp)}{\N\fp}\Bigr) < \pr{\fp \in \cI_r} \Bigl(1 + \frac{k^2(f(\fp, J) - f(\fp, j))^2}{\N(\fp)}  \Bigr) \\
&< \exp \biggl( k^2 \su{\fp \in \cI_r} \frac 1 {\N\fp} \min \Big\{ 1, \frac{4\log^2 \N\fp}{(\theta_j \log X)^2} \Big\} \biggr).
}

\subsubsection{Estimate of $\tilde \cE_2(r,c)$} 
Using Remark \ref{rem:nu}, the fact that \eqs{\max \{\Omega (\fa_{1r}\fb_{1r}), \Omega (\fa_{2r}\fb_{2r}), \Omega (\fa_{3r}), \Omega(\fb_{3r}) \} >  \ell_r}
and the inequality $f(\fa,\cdot) \le 1$, we see that
\bs{
|\tilde \cE_2 (r,c)| 
&< \frac 2 {2^{\ell_r}}
\mathop{\sum \ldots \sum}_{\substack{
\fa_{1r},\fb_{1r}, \fa_{2r},\\
\fb_{2r}, \fa_{3r},\fb_{3r} \in \cI_r  \\ 
\fm_r = \cube }} 
\frac{k^{\Omega (\fa_{1r} \fb_{1r} \fa_{2r}  \fb_{2r})}  2^{\Omega (\fa_{3r}  \fb_{3r})} }
{\N(\fa_{1r}\fb_{1r}\fa_{2r}\fb_{2r} \fa_{3r}\fb_{3r})^{1/2}\kappa^{\Omega(\fa_{3r}\fb_{3r})}}  
\\
& \qquad \cdot \nu(\fa_{1r})\nu(\fb_{1r}) \nu(\fa_{2r})\nu(\fb_{2r})
\nu_{k\kappa} (\fa_{3r}) \nu_{k\kappa} (\fb_{3r}) \\
&= \frac 2 {2^{\ell_r}} \mathop{\sum \ldots \sum}_{\substack{ \fa,\fb, \fa_{3r},\fb_{3r} \in \cI_r  \\  		\fa\fb^2\fa_{3r}\fb_{3r}^2  = \cube }} 
\frac{ k^{\Omega (\fa \fb )} 2^{\Omega (\fa_{3r}  \fb_{3r})}\nu_2 (\fa)\nu_2(\fb) \nu_{k\kappa} (\fa_{3r}) \nu_{k\kappa} (\fb_{3r})} {\N(\fa \fb \fa_{3r}\fb_{3r})^{1/2} \kappa^{\Omega(\fa_{3r}\fb_{3r})}}     \\
&= \frac 2 {2^{\ell_r}} \su{\fA, \fB \in \cI_r \\ \fA\fB^2 = \cube}
\frac{ H(\fA) H(\fB)   }
{\N(\fA \fB)^{1/2}}  = \frac 2 {2^{\ell_r}}\su{\fC \in \cI_r} \frac 1 
{\N(\fC)^{1/2}} \su{\fA, \fB \in \cI_r \\ \fA\fB = \fC \\ \fA\fB^2 = \cube} H(\fA) H(\fB) .
}
where 
\bs{H(\fA) &= \su{\fa, \fa_{3r} \\ \fa\fa_{3r} = \fA} \nu_2 (\fa)  k^{\Omega (\fa)}  \nu_{k\kappa} (\fa_{3r}) (2/\kappa)^{\Omega(\fa_{3r})} .
}
This gives
\bs{
2^{\ell_r-1} |\tilde \cE_2 (r,c)|  &< \prod_{\fp \in \cI_r}  
\biggl(
1+ \sum_{n\ge 2} \frac 1 {\N(\fp)^{n/2}}  \su{0 \le m \le n\\ 3 \mid m+n} H(\fp^{n-m}) H(\fp^m)
\biggr) \\
&< \prod_{\fp \in \cI_r}  \biggl(
1+ \sum_{n\ge 2} \frac{ (8k)^n }
{n! \N(\fp)^{n/2}}  \biggr)  }
since
\eqs{
H(\fp^m) = \su{0 \le i \le m} \nu_2 (\fp^i)  k^{\Omega (\fp^i)}  \nu_{k\kappa} (\fp^{m-i}) (2/\kappa)^{\Omega(\fp^{m-i})} = \frac{(4k)^m}{m!}.
}
Therefore, since $\N\fp > k_0 \ge (2k)^2$, 
\bsc{\label{E2(r,c)bound}
|\tilde \cE_2 (r,c)| <  2^{1-\ell_r} \exp \biggl( 4k^2( e^4 - 5 )\sum_{\fp \in \cI_r} \frac 1 {\N\fp} \biggr). }

Combining \eqref{E2(r,c)bound} with \eqref{E1(r,c)bound} and using Lemma \ref{primesumestimates} it follows that for sufficiently large $X$, 
\eqn{\label{E(r)Productbound}
\prod_{r=0}^j |\tilde \cE (r,c)| <  \exp (8k^2) \prod_{r=0}^j ( 1 + 2^{-\ell_r/2} ). }

\subsection{Bound for $\tilde \cE(j+1)$}
Recall that for $r=j+1$, 
\bs {\tilde \cE(r) 
&= \frac {(2s_r)!} {4^{s_r}} 
\su{\fa_{2r}, \fb_{2r} \in \cI_r\\ 
0\le \Omega (\fa_{2r}), \Omega(\fb_{2r}) \le k\kappa\ell_r
} 
\frac {\lambda(\fa_{2r}\fb_{2r}) f(\fa_{2r}\fb_{2r},J)\nu_{k\kappa} (\fa_{2r}; \ell_r) \nu_{k\kappa} (\fb_{2r}; \ell_r) } {\N(\fa_{2r}\fb_{2r} )^{1/2} \kappa^{\Omega(\fa_{2r}\fb_{2r})}} \\
&\qquad \cdot  \su{\fa_{1r}, \fb_{1r}\in \cI_r 
\\ \Omega (\fa_{1r}\fb_{1r}) = 2s_r \\
\fa_{1r}\fb_{1r}^2\fa_{2r}\fb_{2r}^2 = \cube\\ (c, \prod_i \fa_{ir}\fb_{ir}) = 1}
\frac {f(\fa_{1r}\fb_{1r}, u) \nu(\fa_{1r})\nu(\fb_{1r}) }{\N(\fa_{1r}\fb_{1r})^{1/2}}.
}
We first replace $\fa_{ir}, \fb_{ir}$ with $\fx \fa_{1r}, \fx\fb_{1r}$ and $\fX\fa_{2r}, \fX\fb_{2r}$, respectively, with $(\fa_{ir}, \fb_{ir})=1$ for $i=1,2$, use the inequalities $|\lambda|, f \le 1$, $\nu_n (\fa;\cdot) \le \nu_n (\fa)$, and then remove the condition that $(c, \prod_i \fa_{ir}\fb_{ir}) = 1$. This shows that $4^{s_r} |\tilde \cE(r) | /((2s_r)!)$ is bounded by
\eqs{
\su{\fa_{2r}, \fb_{2r}, \fX \in \cI_r\\ 
(\fa_{2r}, \fb_{2r}) = 1\\
0\le \Omega (\fX\fa_{2r}), \Omega(\fX\fb_{2r}) \le k\kappa\ell_r
} 
\frac {\nu_{k\kappa} (\fX\fa_{2r}) \nu_{k\kappa} (\fX\fb_{2r}) } {\N(\fX^2\fa_{2r}\fb_{2r})^{1/2} \kappa^{\Omega(\fX^2\fa_{2r}\fb_{2r})}} 
\su{\fa_{1r}, \fb_{1r}, \fx \in \cI_r \\ 
(\fa_{1r}, \fb_{1r}) = 1\\
\Omega (\fx^2\fa_{1r}\fb_{1r}) = 2s_r \\
\fa_{1r}\fb_{1r}^2\fa_{2r}\fb_{2r}^2 = \cube}
\frac {\nu(\fx\fa_{1r})\nu(\fx\fb_{1r}) }{\N(\fx^2\fa_{1r}\fb_{1r})^{1/2}}.}
Next, replacing $\fa_{1r}, \fb_{2r}$ with $\fy\fa_{1r}, \fy\fb_{2r}$, and  $\fa_{2r}, \fb_{1r}$ with $\fY\fa_{2r}, \fY\fb_{1r}$, where $(\fa_{ir}, \fb_{3-i,r})=1$ for $i=1,2$ gives
\bs{
|\tilde \cE(r)|	& \le \frac {(2s_r)!} {4^{s_r}} 
\su{\fa_{2r}, \fb_{2r}, \fX, \fy, \fY \in \cI_r\\ 
(\fY\fa_{2r}, \fy\fb_{2r}) = 1\\
0\le \Omega (\fX\fY\fa_{2r}), \Omega(\fy\fX\fb_{2r}) \le k\kappa\ell_r
} 
\frac {\nu_{k\kappa} (\fX\fY\fa_{2r}) \nu_{k\kappa} (\fX\fy\fb_{2r}) } {\N(\fX^2\fy^2\fY^2\fa_{2r}\fb_{2r})^{1/2} \kappa^{\Omega(\fX^2\fy\fY\fa_{2r}\fb_{2r})}} \\
&\qquad \cdot  
\su{\fa_{1r}, \fb_{1r}, \fx \in \cI_r \\ 
(\fy\fa_{1r}, \fY\fb_{1r}) = 1\\
(\fa_{2r}, \fb_{1r})=1 = (\fa_{1r}, \fb_{2r}) \\
\Omega (\fx^2\fa_{1r}\fb_{1r}) = 2s_r - \Omega (\fy\fY) \\
\fa_{1r}\fb_{1r}^2 \fa_{2r}\fb_{2r}^2 = \cube}
\frac {\nu(\fx\fy\fa_{1r})\nu(\fx\fY\fb_{1r}) }{\N(\fx^2\fa_{1r}\fb_{1r})^{1/2}}.}
Now, replace $\fa_{1r}, \fa_{2r}$ with $\fz\fa_{1r}, \fz\fa_{2r}$ and 
$\fb_{1r}, \fb_{2r}$ with $\fZ\fb_{1r}, \fZ\fb_{2r}$,  where $(\fa_{1r}, \fa_{2r})=1= (\fb_{1r}, \fb_{2r})$ to get
\bs{
|\tilde \cE(r)|	& \le 
\frac {(2s_r)!} {4^{s_r}} 
\su{\fa_{2r}, \fb_{2r}, \fX, \fy, \fY,\fz, \fZ \in \cI_r\\ 
(\fY\fz\fa_{2r}, \fy\fZ \fb_{2r}) = 1\\
0\le \Omega (\fX\fY\fz\fa_{2r}), \Omega(\fX\fy\fZ \fb_{2r}) \le k\kappa\ell_r
} 
\frac {\nu_{k\kappa} (\fX\fY\fz\fa_{2r}) \nu_{k\kappa} (\fX\fy\fZ \fb_{2r}) } {\N(\fX^2\fy^2\fY^2\fz\fZ \fa_{2r}\fb_{2r})^{1/2} \kappa^{\Omega(\fX^2\fy\fY\fz\fZ\fa_{2r} \fb_{2r})}} \\
&\qquad \cdot  
\su{\fa_{1r}, \fb_{1r}, \fx \in \cI_r \\
(\fa_{1r}, \fa_{2r})=1= (\fb_{1r}, \fb_{2r})\\ 
(\fy\fz\fa_{1r}, \fY\fZ\fb_{1r}) = 1\\
(\fz\fa_{2r}, \fZ\fb_{1r})=1 = (\fz\fa_{1r}, \fZ \fb_{2r}) \\
\Omega (\fx^2\fa_{1r}\fb_{1r}) = 2s_r - \Omega (\fy\fY\fz\fZ) \\
\fz^2\fa_{1r}\fa_{2r},  \fZ \fb_{1r}^2  \fb_{2r}^2 = \cube}
\frac {\nu(\fx\fy\fz\fa_{1r})\nu(\fx\fY\fZ\fb_{1r}) }{\N(\fx^2\fz\fZ\fa_{1r}\fb_{1r})^{1/2}}.}
Finally writing $\fz = \fa\fb^2\fc^3$ and $\fZ = \fA\fB^2\fC^3$, where $(\fa, \fb) = 1 = (\fA, \fB)$, all square-free, shows that $\fz^2\fa_{1r}\fa_{2r} = \fa^2\fb^4\fc^6\fa_{1r}\fa_{2r} = \cube$ if and only if  $\fa^2\fb\fa_{1r}\fa_{2r} = \cube$ iff 
$\fa_{1r} = \fa_1\fb_1^2\fc_1^3$ and $\fa_{2r} = \fa_2\fb_2^2\fc_2^3$, where $\fa = \fa_1 \fa_2$ and $\fb = \fb_1 \fb_2$. Similarly, we get $\fZ\fb_{1r}^2\fb_{2r}^2 = \fA\fB^2\fC^3 \fb_{1r}^2\fb_{2r}^2 = \cube$ iff $\fb_{1r} = \fA_1 \fB_1^2 \fC_1^3$ and $\fb_{2r} = \fA_2 \fB_2^2\fC_2^3$ where $\fA = \fA_1 \fA_2$ and $\fB = \fB_1 \fB_2$. Making these changes leads to the bound
\bsc{\label{Ej+1chgofvars}
&|\tilde \cE(r)|  \le \frac {(2s_r)!} {4^{s_r}} 
\mathop{\sum\ldots \sum}_{\substack{\fc_2, \fC_2,
\fX, \fy, \fY, \fc, \fC \in \cI_r\\ 
} }
\frac {\nu(\fy\fc^3)\nu(\fY\fC^3) \nu_{k\kappa} (\fX\fY\fc^3\fc_2^3) \nu_{k\kappa} (\fX\fy\fC^3\fC_2^3)}
{\N(\fc^6\fC^6\fX^2\fy^2\fY^2\fc_2^3 \fC_2^3)^{1/2} \kappa^{\Omega(\fX^2\fy\fY\fc^3\fC^3 \fc_2^3\fC_2^3)}} \\
&\cdot \su{\fa \in \cI_r \\ \fa \text{ sq-free}\\
\fa = \fa_1 \fa_2} 
\frac {\nu(\fa\fa_1) \nu_{k\kappa} (\fa\fa_2)}
{\N(\fa)^{3/2} \kappa^{\Omega(\fa\fa_2)}}
\su{\fA \in \cI_r \\ \fA \text{ sq-free}\\
\fA = \fA_1 \fA_2} 
\frac {\nu(\fA\fA_1) \nu_{k\kappa} (\fA\fA_2)}{\N(\fA)^{3/2} \kappa^{\Omega(\fA\fA_2)}} 
\su{\fb \in \cI_r \\ \fb \text{ sq-free}\\
\fb = \fb_1 \fb_2 
} 
\frac {\nu(\fb^2\fb_1^2) \nu_{k\kappa} (\fb^2\fb_2^2)}
{\N(\fb)^3 \kappa^{\Omega(\fb^2\fb_2^2)}} \\
&\cdot 
\su{\fB \in \cI_r \\ \fB \text{ sq-free}, \fB = \fB_1 \fB_2\\ 
0 \le \Omega (\fX\fY\fa\fb^2\fc^3\fa_2\fb_2^2\fc_2^3) \le k\kappa\ell_r \\
0 \le \Omega(\fX\fy\fA\fB^2\fC^3 \fA_2 \fB_2^2\fC_2^3)  \le k\kappa\ell_r
} \frac {\nu(\fB^2\fB_1^2) \nu_{k\kappa} (\fB^2\fB_2^2)}
{\N(\fB)^3 \kappa^{\Omega(\fB^2\fB_2^2)}}
\su{\fc_1, \fC_1, \fx \in \cI_r \\
\Omega (\fx^2\fc_1^3\fC_1^3) = N}
\frac {\nu(\fx\fc_1^3)\nu(\fx\fC_1^3
) }{\N(\fx^2\fc_1^3\fC_1^3)^{1/2}},
}
where 
\eqs{N = 2s_r - \Omega (\fy\fY\fa\fb^2\fc^3\fA\fB^2\fC^3\fa_1\fb_1^2\fA_1 \fB_1^2).}
Since 
\eqs{\Omega (\fy\fY\fa\fb^2\fc^3\fA\fB^2\fC^3\fa_1\fb_1^2\fA_1 \fB_1^2) \le \Omega (\fy\fY(\fa\fb^2)^2 (\fA\fB^2)^2\fc^3\fC^3) \le 4k\kappa\ell_r}
we see that $N \ge 2s_r - 4k\kappa\ell_r$. Writing $N = 3q + R$ with $0 \le R < 3$, it follows that 
\bs{
&\su{\fc_1, \fC_1, \fx \in \cI_r \\
\Omega (\fx^2\fc_1^3\fC_1^3) = N 
} 
\frac {\nu(\fx\fc_1^3)\nu(\fx\fC_1^3
) }{\N(\fx^2\fc_1^3\fC_1^3)^{1/2}} \\
&\le \su{a,b,c \ge 0\\ 2a+3b+3c=N}
\su{\fx \in \cI_r \\
\Omega (\fx) = a} 
\frac {\nu(\fx)}{\N(\fx)} 
\su{\fc_1 \in \cI_r \\
\Omega(\fc_1) = b} 
\frac {\nu(\fc_1)}{3^b\N(\fc_1)^{3/2}} 
\su{\fC_1\in \cI_r \\
\Omega(\fC_1) = c} 
\frac {\nu(\fC_1) }{3^c \N(\fC_1)^{3/2}} \\
&= \su{a,b,c \ge 0\\ 2a+3b+3c=N}
\frac 1 {a! b!c!} 
\biggl(\su{\fp \in \cI_r} \frac 1 {\N\fp}\biggr)^{a} 
\biggl(\su{\fp \in \cI_r} \frac 1 {3(\N\fp)^{3/2}}\biggr)^{b+c}.
}
Using Lemma \ref{primesumestimates} we have
\bs{
\su{\fc_1, \fC_1, \fx \in \cI_r \\
\Omega (\fx^2\fc_1^3\fC_1^3) = N 
} 
\frac {\nu(\fx\fc_1^3)\nu(\fx\fC_1^3
) }{\N(\fx^2\fc_1^3\fC_1^3)^{1/2}}  &< \su{0 \le 2a \le 3q+R \\ 3 \mid R - 2a} \frac 1 {a!( q -\frac{2a-R}{3})!} \biggl(\su{\fp \in \cI_r} \frac 1 {\N\fp}\biggr)^{a} 
\Bigl( \frac 2 3 \Bigr)^{(N-2a)/3} \\
&= \su{0 \le m \le q\\ 2 \mid 3m + R} \frac {(2/3)^{q-m}} {(\frac{3m + R}{2})! ( q - m )!} \biggl(\su{\fp \in \cI_r} \frac 1 {\N\fp}\biggr)^{(3m+R)/2} 
\\
&\le \biggl(\su{\fp \in \cI_r} \frac 1 {\N\fp}\biggr)^{s_r}  \frac 1 { q!} \su{0 \le m \le q\\ 2 \mid 3m + R} \binom{q}{m}  
\Bigl(\frac 23 \Bigr)^{q-m}\\
&\le 
\frac 1 {\fl{\frac{2s_r - 4k\kappa\ell_r}{3}}!} \biggl(\su{\fp \in \cI_r} \frac 1 {N\fp}  \biggr)^{s_r}  \Bigl(\frac 5 3\Bigr) ^{\fl{\frac{2s_r - 4k\kappa\ell_r}{3}}}.
}
where in the last inequality we used the fact that $(5/3)^q/q!$ is decreasing for $q\ge 1$ and in our case $q = \fl{N/3} \ge 1$ by the choice of $s_r$. Since the last result is independent of $N$, we can bound the remaining terms in \eqref{Ej+1chgofvars} independent of the innermost sum. The sum over $\fB$ is bounded by
\bs{
&\prod_{\fp \in \cI_r} \biggl( 1 + 
\frac {\nu(\fp^2)\nu_{k\kappa} (\fp^2)} {\N(\fp)^3 \kappa^{\Omega(\fp^2)}} \su{\fB_1, \fB_2\\ \fp = \fB_1\fB_2}
\frac{\nu(\fB_1^2) \nu_{k\kappa} (\fB_2^2) }{\kappa^{\Omega(\fB_2^2)}} 
\biggr) = \prod_{\fp \in \cI_r} \biggl( 1 + 
\frac{k^2  + k^4}{8(\N\fp)^3} 
\biggr).
}
The same bound holds for the sum over $\fb$. As for the sum over $\fa$, similarly for $\fA$, we have 
\bs{
\su{\fa \in \cI_r \\ \fa \text{ sq-free}\\	\fa = \fa_1 \fa_2} \frac {\nu(\fa\fa_1) \nu_{k\kappa} (\fa\fa_2)} {\N(\fa)^{3/2} \kappa^{\Omega(\fa\fa_2)}} 
&\le \su{\fa \in \cI_r \\ \fa \text{ sq-free}\\ } \frac {\nu(\fa) \nu_{k\kappa} (\fa)} {\N(\fa)^{3/2}\kappa^{\Omega(\fa)}}  \su{\fa_1, \fa_2\\ \fa = \fa_1 \fa_2} \frac {\nu(\fa_1)\nu_{k\kappa} (\fa_2)} {\kappa^{\Omega(\fa_2)}} \\
& = \prod_{\fp \in \cI_r} \biggl( 1 + 
\frac{k+k^2}{(\N\fp)^{3/2}}
\biggr) < 
\exp  \biggl( \sum_{\fp \in \cI_r} \frac {k(1+k)} {\N(\fp)^{3/2}} \biggr).
}
Finally we see that 
\bs{
&
\mathop{\sum\ldots \sum}_{\substack{\fc_2, \fC_2,
\fX, \fy, \fY, \fc, \fC \in \cI_r\\ 
} }
\frac {\nu(\fy\fc^3)\nu(\fY\fC^3) \nu_{k\kappa} (\fX\fY\fc^3\fc_2^3) \nu_{k\kappa} (\fX\fy\fC^3\fC_2^3)}
{\N(\fc^6\fC^6\fX^2\fy^2\fY^2\fc_2^3 \fC_2^3)^{1/2} \kappa^{\Omega(\fX^2\fy\fY\fc^3\fC^3 \fc_2^3\fC_2^3)}} \\
&\le
\biggl(\sum_{\fy \in \cI_r}
\frac {\nu(\fy) \nu_{k\kappa} (\fy) } {\N(\fy) \kappa^{\Omega(\fy)}}
\sum_{\fc_2  \in \cI_r}
\frac {\nu_{k\kappa} (\fc_2^3) } {\N(\fc_2)^{3/2} \kappa^{\Omega(\fc_2^3)}}
\sum_{\fc  \in \cI_r}
\frac {\nu(\fc^3)\nu_{k\kappa} (\fc^3) }{\N(\fc)^3 \kappa^{\Omega(\fc^3)}}\biggr)^2
\sum_{\fX  \in \cI_r}
\frac {\nu_{k\kappa} (\fX) \nu_{k\kappa} (\fX) } {\N(\fX) \kappa^{\Omega(\fX^2)}}.
}
The sum over $\fX$ is
\bs{
& 
\sum_{\fX  \in \cI_r}
\frac {\nu_{k\kappa} (\fX) \nu_{k\kappa} (\fX) } {\N(\fX) \kappa^{\Omega(\fX^2)}} = \prod_{\fp \in \cI_r} \Bigl( 1 + \sum_{n\ge 1} 
\frac {\nu_{k\kappa} (\fp^n) \nu_{k\kappa} (\fp^n) } {\N(\fp^n) \kappa^{\Omega(\fp^{2n})}} 
\Bigr) \\
&= \prod_{\fp \in \cI_r} \Bigl( 1 + \sum_{n\ge 1} 
\frac { k^{2n}  } {(n!)^2 \N(\fp)^n} 
\Bigr) < \exp \biggl( e^{1/4}k^2 \sum_{\fp \in \cI_r}  \frac {1} {\N(\fp)}\biggr)
}
since $k^2/\N\fp < 1/4$ and for $0 \le x \le 1/4$, 
\eqs{\sum_{n \ge 1} \frac{x^n}{(n!)^2} \le e^{1/4} x.}
Similarly, we see that 
\bs{\sum_{\fc_2  \in \cI_r}
\frac {\nu_{k\kappa} (\fc_2^3) } {\N(\fc_2)^{3/2} \kappa^{\Omega(\fc_2^3)}}  &= \prod_{\fp \in \cI_r} \Bigl( 1 + \sum_{n\ge 1} 	\frac {k^{3n}  } {(3n)! \N(\fp)^{3n/2} } 	\Bigr) < \exp \biggl(  \sum_{\fp \in \cI_r} \frac {0.003k^3} {\N(\fp)^{3/2}}\biggr) \\
\sum_{\fy \in \cI_r}
\frac {\nu(\fy) \nu_{k\kappa} (\fy) } {\N(\fy) \kappa^{\Omega(\fy)}} 
&= \prod_{\fp \in \cI_r} \Bigl( 1 + \sum_{n\ge 1}  \frac {k^n  } {(n!)^2 (\N\fp)^n }  \Bigr) < \exp \biggl( \sum_{\fp \in \cI_r} \frac{ke^{1/(4k)}} {\N\fp}\biggr)\\
\sum_{\fc  \in \cI_r}
\frac {\nu(\fc^3)\nu_{k\kappa} (\fc^3) }{\N(\fc)^3 \kappa^{\Omega(\fc^3)}} 
&= \prod_{\fp \in \cI_r} \Bigl( 1 + \sum_{n\ge 1} \frac {k^{3n}  } {((3n)!)^2 \N(\fp)^{3n} } \Bigr) < \exp \biggl(  \sum_{\fp \in \cI_r} \frac{0.003k^3} {\N(\fp)^3}\biggr),
}
by using the inequality
\eqs{\sum_{n\ge 1} 
\frac {x^n } {(3n)!} < x \sum_{n\ge 1} \frac 1 {(3n)! 8^{n-1}} \qquad (0 < x \le 1/8)
}
and the fact that since $\omega^2 + \omega + 1 = 0$, 
\bs{8\sum_{n\ge 1} \frac 1 {(3n)!8^n} &= \frac 8 3 \( \sum_{n\ge 1}  \frac{(\omega/8)^n}{n!} + \sum_{n\ge 1}  \frac{(\omega^2/8)^n}{n!} + \sum_{n\ge 1}  \frac 1 {n! 8^n} \)  \\
&= \frac 8 3 (e^{1/8} + e^{\omega/8} + e^{\omega^2/8} - 3)  \approx 0.002604.}
Combining all these estimates and using Lemma \ref{primesumestimates} we conclude that 
\bsc{\label{E(r,c)bound2}
&|\tilde \cE(r,c)| < 
\frac {(2s_r)!} {4^{s_r}\fl{\frac{2s_r - 4k\kappa\ell_r}{3}}!} \biggl(\su{\fp \in \cI_r} \frac 1 {N\fp}  \biggr)^{s_r}  \Bigl(\frac 5 3\Bigr)^{\fl{\frac{2s_r - 4k\kappa\ell_r}{3}}} \\
&\cdot 
\exp \biggl( \frac{0.024k^3 + k^2  + k^4}{4\fl{X^{\theta_{r-1}}}^2}   
+ \sum_{\fp \in \cI_r}  \frac{0.006k^3 + 2k+ 2k^2 }{(\N\fp)^{3/2}}
+    \frac {e^{1/4}k^2 + 2ke^{1/(4k)}} {\N\fp}
\biggr).
}

\subsection{Bound for $\tilde \cE(r,c)$ for $j+1 < r \le J$}
Recall that for $j+1 < r \le J$,
\bs{\tilde \cE(r,c) &= \su{\fa_{1r}, \fb_{1r} \in \cI_r\\\fa_{1r}\fb_{1r}^2 = \cube \\ (c,\fa_{1r}\fb_{1r})=1} \frac{\lambda(\fa_{1r}\fb_{1r}) f(\fa_{1r}\fb_{1r},J)}
{\N(\fa_{1r}\fb_{1r})^{1/2}\kappa^{\Omega(\fa_{1r}\fb_{1r})}}
\nu_{k\kappa} (\fa_{1r}; \ell_r) \nu_{k\kappa} (\fb_{1r}; \ell_r).}
Taking absolute values we have as before
\bsc{\label{E(r,c)bound3}
|\tilde \cE(r,c)| &\le 
\su{\fa_{1r}, \fb_{1r} \in \cI_r\\ \fa_{1r}\fb_{1r}^2 = \cube} 
\frac{ \nu_{k\kappa} (\fa_{1r}) \nu_{k\kappa} (\fb_{1r})}
{\N(\fa_{1r}\fb_{1r})^{1/2}\kappa^{\Omega(\fa_{1r}\fb_{1r})}}  \\
&= \su{\fa \in \cI_r} 
\frac {1}{\N(\fa)^{1/2}\kappa^{\Omega(\fa)}}  
\su{\fa_{1r}, \fb_{1r} \in \cI_r \\ \fa = \fa_{1r}\fb_{1r} \\ \fa\fb_{1r} = \cube } \nu_{k\kappa} (\fa_{1r}) \nu_{k\kappa} (\fb_{1r}) \\
&= \pr{\fp \in \cI_r} \biggl(1 + 
\sum_{n \ge 1} 
\frac{1}{\N(\fp)^{n/2} \kappa^n}
\su{0 \le m \le n \\ 3 \mid m+n} \nu_{k\kappa} (\fp^{n-m}) \nu_{k\kappa} (\fp^m)
\biggr) \\
&< \pr{\fp \in \cI_r} \biggl(1 + 
\sum_{n \ge 2} 
\frac{(2k)^n}{n! \N(\fp)^{n/2}} \biggr) \le \exp \biggl( 4(e-2)k^2 \su{\fp \in \cI_r}  \frac{1} {\N\fp}  \biggr).
}

\subsection{Moments of $L(s,\chi_c)$} 

\begin{lem} \label{2lthmoment}
For any complex numbers $z_\fp$ with $|z_\fp| \le 1$, any positive integer $l = 3q + r$ with $0 \le r < 3$, and any $\eps > 0$, 
\bs{
\sum_{c \in \sF(X)} \bigg| \sum_{3< \N\fp \le y} \frac{\chi_c (\fp)z_\fp }{\N(\fp)^{1/2}}   \bigg|^{2l} &< |\sF(X)|
\biggl(\max \biggl\{ 1, \su{3 < \N\fp \le y} \frac{|z_\fp|^2}{\N\fp}  \biggr\} \biggr)^l \frac{(l!)^2( 19/9)^{2(l-r)/3}}{(2(l-r)/3)!}  \\
&\qquad + O \biggl(X^{1/2+\eps} \biggl( \sum_{3< \N\fp \le y} \frac{|z_\fp|}{\N(\fp)^{1/2-\eps}} \biggr)^{2l}\biggr).
}
\end{lem}
\begin{proof}
Put $\cI = (3,y]$. Then, 
\bs{ \bigg| \sum_{\fp \in \cI} \frac{\chi_c (\fp) z_\fp}{\N(\fp)^{1/2}}    \bigg|^{2l} 
&= (l!)^2 \su{\fa_1, \fa_2 \in \cI\\ \Omega(\fa_i) = l } \frac{\chi_c (\fa_1 \fa_2^2)}{\N(\fa_1\fa_2)^{1/2}} \nu (\fa_1) \nu (\fa_2) z_{\fa_1} \overline{z_{a_2}},
}
and summing over $c \in \sF(X)$ we get
\eqs{(l!)^2 \su{\fa_1, \fa_2 \in \cI\\ \Omega(\fa_i) = l } \frac{\nu (\fa_1) \nu (\fa_2) z_{\fa_1} \overline{z_{a_2}} }{\N(\fa_1\fa_2)^{1/2}} 
\su{c \in \sF(X) \\ (c,\fa_1\fa_2)=1}  \chi_c (\fa_1\fa_2^2).
}
When $\fa_1\fa_2^2 =\cube$, this sum is trivially bounded by
\eqs{|\sF(X)| (l!)^2 \su{\fa_1, \fa_2 \in \cI\\ \Omega(\fa_i) = l \\ \fa_1\fa_2^2 =\cube} \frac{\nu (\fa_1) \nu (\fa_2) |z_{\fa_1}| |z_{a_2}|}{\N(\fa_1\fa_2)^{1/2}}.
}
Otherwise, applying Lemma \ref{charsumoverfamily} gives 
\bs{\su{c \in \sF(X) \\ (c,\fa_1\fa_2)=1}  \chi_c (\fa_1\fa_2^2) \ll X^{1/2+\eps} \N(\fa_1\fa_2)^\eps,}
which then leads to the claimed error term. For the case of cubes, replacing $\fa_1, \fa_2$ by $\fb\fa_1, \fb\fa_2$ with $(\fa_1, \fa_2)=1$ we see that $\fa_1\fb \fa_2^2 \fb^2 =\cube$ if and only if $\fa_1, \fa_2$ are cubes. Thus, 
\bs{
&\su{\fa_1, \fa_2 \in \cI\\ \Omega(\fa_i) = l \\ \fa_1\fa_2^2 =\cube} \frac{\nu (\fa_1) \nu (\fa_2) |z_{\fa_1}| |z_{a_2}|}{\N(\fa_1\fa_2)^{1/2}} 
\le \su{\fb \in \cI\\ 0 \le \Omega(\fb) \le l\\ 3 \mid l - \Omega(\fb)} \frac {|z_{\fb}|^2 \nu(\fb)} {\N(\fb)} 
\biggl(\su{\fa \in \cI\\ \Omega(\fb\fa^3) = l} \frac{\nu (\fa^3) |z_{\fa}|^3  }{\N(\fa)^{3/2}} \biggr)^2 \\
&< \su{\fb \in \cI\\ 0 \le \Omega(\fb) \le l\\ 3 \mid l - \Omega(\fb)} \frac {|z_{\fb}|^2 \nu(\fb)} { 9^{(l-\Omega(\fb))/3}  \bigl((\frac{l-\Omega(\fb)} 3)!\bigr)^2 \N(\fb)} 
}
since
\bs{
\su{\fa \in \cI\\ \Omega(\fb\fa^3) = l} \frac{\nu (\fa^3) |z_{\fa}|^3  }{\N(\fa)^{3/2}}  &\le 3^{-(l-\Omega(\fb))/3} \su{\fa \in \cI\\ \Omega(\fa) = (l-\Omega(\fb))/3} \frac{\nu(\fa) }{\N(\fa)^{3/2}}  \\
& = \frac 1 {3^{(l-\Omega(\fb))/3} ( \frac{l-\Omega(\fb)} 3 )!} \biggl(\su{\fp \in \cI} \frac 1 {\N(\fp)^{3/2}}  \biggr)^{(l-\Omega(\fb))/3}
}
and by Lemma \ref{primesumestimates} the last sum over primes is less than 1. Here, we also used the inequalities $\nu(ab) \le \nu(a)\nu(b)$, $\nu(b)^2 \le \nu(b)$, $\nu(a^3) \le 3^{-\Omega(a)} \nu(a)$, and ignored the condition that $(a_1,a_2) = 1$. Writing $l = 3q + r$ with $0\le r < 3$, the sum over $b$ becomes 
\bs{
& \su{0 \le j \le q}  
\frac 1 {9^{q-j} \bigl(\bigl(q - j\bigr)!\bigr)^2 }  \su{\fb \in \cI\\ \Omega(b) = 3j+r} \frac{|z_\fb|^2 \nu (\fb)}{\N(\fb)}\\
&= \su{0 \le j \le q}  
\frac 1 {9^{q-j} \bigl(\bigl(q - j\bigr)!\bigr)^2 (3j+r)!}  \biggl(\su{3 < \N\fp \le y} \frac{|z_\fp|^2}{\N\fp}\biggr)^{3j+r} 
\\
& \le \biggl(\max \biggl\{ 1, \su{3 < \N\fp \le y} \frac{|z_\fp|^2}{\N\fp}  \biggr\} \biggr)^l \frac{( 1+ 1+ 1/9)^{2(l-r)/3}}{(2(l-r)/3)!}.
}
This gives the desired result.
\end{proof}

The next result used in the proof of Theorem \ref{thm:upperbound} is not optimal, but is sufficient for our purposes in this paper. 
\begin{prop} \label{prop:2kthmoment}
For a real number $k>0$,
\eqs{\su{c \in \sF(X)} |L(1/2,\chi_c)|^{2k} \ll |\sF(X)| \exp\( 2k (\log \log X)^2 \).}
\end{prop}
\begin{proof}
Put 
\eqs{N(v) = | \{ c \in \sF(X) : L(1/2,\chi_c) \neq 0 \wedge \log |L(1/2,\chi_c)| \ge v \}|.
}
Then, since trivially $N(v) \le |\sF(X)|$,
\bsc{\label{2kmomentintegral}
\su{c \in \sF(X)} |L(1/2,\chi_c)|^{2k} &= \int_\infty^{-\infty} e^{2kv} dN(v) = 2k \int_{-\infty}^\infty e^{2kv} N(v) dv 
\\ 
&= 2k \int_z^\infty e^{2kv} N(v) dv + O \( e^{2kz} |\sF(X)| \) 
}
where $z = (\log \log X)^2$. 

Taking $\lambda = 1$ and $t=0$ in Lemma \ref{lem:logderLbound} shows that for any $x \ge 3$, 
\bsc{\label{spongebob}
\log |L(1/2,\chi_c)| &< \re \sum_{1 < \N\fa \le x} \frac{\Lambda(\fa) \chi_c (\fa)}{\N(\fa)^{1/2+1/\log x}\log \N\fa}  \frac{\log (x/\N\fa)}{\log x} \\
&\quad + \frac {\log \cond{c}} {\log x} + \frac 3 {x^{1/2}\log^2 x}.}
If $L(1/2,\chi_c) \neq 0$ and $\log |L(1/2,\chi_c)| \ge v$ for some $c \in\sF(X)$, then taking $x=3$ yields
\eqs{v \le \log |L(1/2,\chi_c)| < \frac {\log \cond{c}} {\log 3} + \frac {3^{1/2}} {\log^2 3} < 1.44 + \log X.}
Thus, we shall assume from now on that $z \le v < \log X + 2$. 

Choose $x $  such that 
\eqs{\log x = \frac A v\log X, \qquad A = A(v) = \log \log v.}
Note that 
\eqs{\log x > \frac{ \log \log (2+\log X)}{2+\log X } \log X \gg \log \log \log X}
as $A(v) / v$ is decreasing for $v \ge z$. 

Since 
\eqs{\su{(\N\fp)^k \le x\\k \ge 3} \frac 1{k\N\fp^{k/2}} \le \frac 1 {3( 1- 3^{-1/2})} \sum_{\N\fp \ge 3} \frac 1 {(\N\fp)^{3/2}} < 1,}
it follows from \eqref{spongebob} that
\bs{v \le \log |L(1/2,\chi)| 
&< \re \sum_{3< \N\fp \le x} \frac{\chi (\fp)}{\N(\fp)^{1/2+1/\log x}}  \frac{\log (x/\N\fp)}{\log x} \\
&\quad + \re \sum_{3< \N\fp \le \sqrt x} \frac{\chi (\fp^2)}{2\N(\fp)^{1+2/\log x}}  \frac{\log (x/\N\fp^2)}{\log x} + \frac v A + O(1).} 
Furthermore, as $\sqrt x \ge 3$ is satisfied with our choice of $x$, we have 
\bs{\sum_{3 < \N\fp \le \sqrt x} \frac{\chi (\fp^2)}{2\N(\fp)^{1+2/\log x}}  \frac{\log (x/\N\fp^2)}{\log x}  
&\le \sum_{p \le \sqrt x} \frac 1 p = \log \log x + O(1), }
which implies that
\bs{
\re \sum_{3 < \N\fp \le x} &\frac{\chi_c (\fp)}{\N(\fp)^{1/2+1/\log x}}  \frac{\log (x/\N\fp)}{\log x}  > v  - \frac v A  - \log \log x + O(1) > v-\frac {2v}A > 0 }
for $v \ge z$ since
\eqs{\log \log x <  \log \log X < \frac z {\log \log z} \le \frac v A }
for sufficiently large $X$ ($X > 181$ works). Applying Lemma \ref{2lthmoment} with \eqs{z_\fp = \frac{\log (x/\N\fp)}{\log x \N(\fp)^{1/\log x}}}
we derive the estimate
\bs{N(v) &< v_1^{-2l} \su{c \in \sF(X)}
\bigg| \sum_{3<\N\fp \le x} \frac{\chi_c (\fp)z_\fp}{\N(\fp)^{1/2}}  \bigg|^{2l} \\
& \ll v_1^{-2l} |\sF(X)| (2\log \log x)^l \frac{(l!)^2( 19/9)^{2(l-r)/3}}{(2(l-r)/3)!}  + v_1^{-2l} X^{1/2+\eps} x^{2l(1/2+\eps)}, }
where $v_1 = v-2v/A$. Using \eqref{Stirlingbounds} gives for $l \ge 3$ and $3 \mid l$ that 
\bs{
\frac{(l!)^2(19/9)^{2(l-r)/3}}{(2(l-r)/3)!} 
&< \( \frac {19} {6e^2}\)^{2l/3} l^{4l/3}  \sqrt{3\pi l} \exp\( \frac 1 {6l} -\frac 1 {8l+1} \) < l^{4l/3}.
}
Thus, choosing $l = 3 \fl{v/(9A)}$, we deduce for $\eps < 1/(12)$ that  
\bsc{\label{Nvbound}
N(v) 
&\ll|\sF(X)| \exp \Bigl( l  \bigl(\log (2\log\log x) + \frac 4 3 \log l - 2 \log v_1 \bigr) \Bigr)  \\	&\quad + X \exp \Bigl( 2l   (1/2+\eps) \log x -  2l \log v_1 - (1/2-\eps) \log X \Bigr) \\
&\ll |\sF(X)| \exp \Bigl(  -\frac {v\log v} {22A}  \Bigr) + X \exp \Bigl(-\frac{v\log v} {2A}\Bigr) }
since using the inequality $-\log(1-x) < x/(1-x)$ we have that
\bs{ \log (2\log\log x) + \frac 4 3  \log l - 2 \log v_1 &< \log 2 + \frac 1 2 \log z  - \frac 2 3 \log v + \frac 4 {A  - 2} \\
&\le -\frac 1 6 \log v + \frac 4 {\log \log z  - 2} < -\frac 1 7 \log v, }
and
\bs{
&2l   (1/2+\eps) \log x -  2l \log v_1 - (1/2-\eps) \log X \\
&<  \frac {2v}{3A}(1/2+\eps) \frac A v \log X -  2l \log v_1 - (1/2-\eps) \log X \\
&< (2\eps- 1 /6 )  \log  X - 6\fl{v/(9A)} \log v_1 < - 6\fl{v/(9A)} \(\log v  - \frac 2 {A-2}\) \\
& < - 5 (v/(9A) - 1) \log v < -\frac{v\log v} {2A}.
}
Inserting \eqref{Nvbound} in \eqref{2kmomentintegral} we get the claimed result. 
\end{proof}

\subsection{Squares of Primes} \label{sec:squares}
In this part we prove the following result.
\begin{lem}
For $0 \le j \le J$ and $X$ satisfying \eqref{Xbound}, the estimate
\eqn{\label{Sjkbound}
\sum_{c \in \sF(X)} (S_{j,k} (c))^2 \ll \cS_k |\sF(X)|  }
holds, where 
\eqn{\label{Sk}
\cS_k =   4\cdot \w^{2kD}(2kD)! \exp \Bigl( 4k(1+D) \Bigr), 
}	
with $D$ given in \eqref{constantD}, provided that 
\eqn{\label{thetacond5}
\red{\Theta  < \frac 1 {8kD + 4\log (16 \sqrt\w)}}.}
\end{lem}
\begin{rem}
For $k=2, \kappa = 1$, one can take $D=1$, in which case
\eqs{\cS_2 \approxeq 5.3316663123\times 10^{11}. }
This follows by using the inequality (see \cite[Section 3, Theorem 5]{RoSc})
\eqs{\sum_{x < p \le 5x} \frac 1 p < \log \log (5x)  - \log \log x + \frac 1 {\log^2 x}, }
which holds for $x \ge 286$ (where we take $x= 8\cdot 5^{n-1}$ and $n \ge 3$), and the fact that $n \sum_{p \in \cI_n} 1/p$ for $n=1, 2$ are $0.416533, 0.353$, respectively. 
\end{rem}
\begin{proof}
Recall that $S_{j,k} (c)  = \exp \bigl(k \re K(c,j) \bigr)$, 
where
\eqs{
K(c,j) = \sum_{\fp \in \cI} \frac{\chi_c(\fp^2) B_j(\fp) }{\N\fp}, \qquad  
B_j(\fp) = \frac 1 {2\N(\fp)^{2/(\theta_j\log X)}}  \Bigl(1 - \frac{2\log \N\fp}{\theta_j\log X}\Bigr),
}	
and $\cI = (k_0, X^{\theta_j/2}]$. Note that $B_j (\fp) \le 1/2$. We extend $B_j$ to be a completely multiplicative function below.

Define the intervals
\eqs{\cI_m = ( k_0 \w^{m-1}, k_0 \w^m ] \cap \cI}
for $m=1, \ldots, M$ with
\eqs{
M = \ceil{\frac{\frac 1 2 \theta_j \log X - \log k_0}{\log \w}}.}
Let
\eqs{
Q_m = \{ c \in \sF(X) : |\re K_m (c,j)| > 2^{-m}, |\re K_n (c,j)| \le 2^{-n} \text{ for } m < n \le M \}}
where
\eqs{K_n (c,j) = \sum_{\fp \in \cI_n} \frac{\chi_c(\fp^2) B_j(\fp) }{\N\fp}. }
Then,
\eqs{
\sum_{c \in \sF(X)} (S_{j,k} (c))^2 \le 
\su{c \in \sF(X) \setminus \cup_m Q_m} (S_{j,k} (c))^2 + \sum_{m=1}^M \su{c \in Q_m} (S_{j,k} (c))^2.
}
For the first sum on the right, we have
\eqn{\label{Squarebound1}
\su{c \in \sF(X) \setminus \cup_m Q_m} (S_{j,k} (c))^2 \le e^{2k}  |\sF(X)| . 
}
If $c \in Q_m$ for some $m$, then we have
\eqs{
\re K(c,j) = \sum_{n=1}^M \re K_n (c,j)  < \sum_{n=1}^m \frac D n + \sum_{n > m} 2^{-n} < D (2 + \log m) + \frac 1 {2^{m}},
}
where we used the estimate
\eqn{\label{constantD}
\sum_{\fp \in \cI_n} \frac{B_j (\fp)}{\N\fp} \le \sum_{p \in \cI_n} \frac 1 p < \frac D n
}
for some positive integer $D$. Therefore, 
\eqs{\su{c \in Q_m} (S_{j,k} (c))^2  \le \su{c \in \sF(X)} \exp \biggl(2k \Bigl( \frac 1 {2^{m}} + D (2 + \log m ) \Bigr)\biggr)  (2^m \re K_m (c,j))^4.}
Noting that
\eqs{
(\re K_m(c,j))^4 = \frac {4!} {2^4} \su{\fa, \fb \in \cI_m \\ \Omega (\fa\fb) = 4 } 
\frac {\nu(\fa) \nu(\fb)  B_j(\fa\fb)}{\N(\fa\fb)} \chi_c (\fa\fb^2),
}
and using Lemma \ref{charsumoverfamily} we deduce that
\bs{
&\sum_{c \in \sF(X)} (\re K_m(c,j))^4
= \frac 32 \su{\fa, \fb \in \cI_m \\ \Omega (\fa\fb) = 4 } 
\frac {\nu(\fa) \nu(\fb)  B_j(\fa\fb)}{\N(\fa\fb)} \sum_{c \in \sF(X)} \chi_c (\fa\fb^2) \\
&\ll  |\sF(X)|
\su{\fa, \fb \in \cI_m \\ \Omega (\fa\fb) = 4 \\ \fa\fb^2 = \cube } 
\frac {\nu(\fa) \nu(\fb)  B_j(\fa\fb)}{\N(\fa\fb)} + X^{1/2+\eps} \biggl(  \sum_{\fp \in \cI_m} \frac{B_j(\fp) }{(\N\fp)^{1-\eps}} \biggr)^4.
}
Now replacing $\fa, \fb$ by $\fc\fa^3, \fc\fb^3$ with $(\fa, \fb)=1$, using $B_j (\fp) \le 1/2$, we see that
\bs{&
\su{\fa, \fb \in \cI_m \\ \Omega (\fa\fb) = 4 \\ \fa\fb^2 = \cube } 
\frac {\nu(\fa) \nu(\fb)  B_j(\fa\fb)}{\N(\fa\fb)} 
\le \su{\fc \in \cI_m \\ \Omega(\fc) = 2} 
\frac{\nu(\fc) B_j (\fc^2)}{(\N\fc)^2} = \frac 1 2 \biggl( \su{\fp\in \cI_m}  \frac{B_j(\fp^2) }{\N(\fp)^2}\biggr)^2 \\
&\le \frac 1 8 \biggl(\sum_{p > \w^{m-1} k_0} \frac 1 {p^2} \biggr)^2
\le \frac 1 8 \biggl( \int_{3\cdot\w^{m-1}}^\infty x^{-2} dx  \biggr)^2 = \frac 1 {72\cdot \w^{2m-2}}.
}
Using this bound together with the estimate
\eqs{
\sum_{\fp \in \cI_m} \frac{B_j(\fp) }{(\N\fp)^{1-\eps}} \le \sum_{p \in \cI_m} \frac 1 {p^{1-\eps}} \le \int_1^{k_0 \w^m} x^{\eps-1} dx < \frac 1 \eps (k_0 \w^m)^\eps, 
}
we conclude that
\mult{\sum_{1 \le m \le M} \su{c \in Q_m} (S_{j,k} (c))^2  \ll |\sF(X)| \sum_{1 \le m \le M} \exp \biggl(2k \Bigl( \frac 1 {2^{m}} + D (2+ \log m ) \Bigr)\biggr)  (4/\w) ^{m}  \\
+ X^{1/2+\eps} \sum_{1 \le m \le M} \exp \biggl(2k \Bigl( \frac 1 {2^{m}} + D (\log m + 2) \Bigr)\biggr)  2^{4m}  \biggl(  \frac 1 \eps (k_0 \w^m)^\eps \biggr)^4.  }
The first term on the right is bounded by
\eqn{
\label{Squarebound2}	
4\cdot \w^{2kD}(2kD)! \exp \Bigl( 2k  + 4kD \Bigr) |\sF(X)|,}
where we used the inequality
\eqs{\sum_{m \ge 1} m^A x^m \le \frac{A! x}{(1-x)^{A+1}}}
which holds for $x < 1$.

The second term is bounded by
\eqs{
\exp \biggl( (1/2+\eps)\log X + 2k \Bigl( 1 + D (2+ \log M) \Bigr) + \log  \frac{ \eps^{-4} k_0^{4\eps} } {2^4 \w^{4\eps} -1} + (M+1) \log (2^4 \w^{4\eps} ) \biggr).
}
We choose $\eps = 1/8$. Then, for $X$ satisfying 
\eqn{\label{Xbound}
2k +4kD + \log  \frac{ 8^4 k_0^{1/2} } {16 \sqrt{\w} -1}  < \frac 18 \log X,
}
the second term is $\ll X$ provided that \eqref{thetacond5} holds. Combining \eqref{Squarebound1} and \eqref{Squarebound2} the claim in \eqref{Sjkbound} follows.
\end{proof}

\section{First Mollified Moment} 

\subsection{Approximate Functional Equation}

\begin{lem} 
\label{lem:approxfnceqn}
Let $\chi$ be a primitive ray class character, $G(u)$ be any even function, holomorphic and bounded in the the region $|\re(u)| < 3$ with $G(0)=1$, and $Y > 0$ be a constant real number. Then, for any complex number $s$ with $\re(s) \in [0,1]$,  
\mult{L(s,\chi) = \sum_{\fb}  \chi (\fb) (\N\fb)^{-s} V_s \biggl( \frac {\N\fb} {Y\sqrt{3\cond{\chi} }} \biggl)  \\
+  (2\pi)^{2s-1} (3\cond{\chi})^{(1/2-s)} \frac{W(\chi)}{\sqrt{\cond{\chi}}}\sum_\fb \frac{\overline{\chi (\fb)}}{(\N\fb)^{1-s}}   V_{1-s} \biggl( \frac {Y\N\fb} {\sqrt{3\cond{\chi}}} \biggr), }
where 
\eqs{\label{Vs(y)}
V_s (y) = \frac 1 {2\pi i} \int_{2-i\infty}^{2+i\infty} G(u) (2 \pi y)^{-u} \frac{\Gamma(s+u)}{\Gamma(s)} \frac{du}{u}.
}
\end{lem}
\begin{rem} 
For $G(u)=1, s=1/2, 0< \alpha \le 1/2-\eps$ with $\eps < 1/2$, and $A > 0$, it follows from Stirling's formula for the Gamma function (see, for example, \cite[Eqn. 5.112]{IwaKow}) by shifting the contour to $\re u = -\alpha$ and to $\re u = A$, respectively, that 
\bsc{\label{Vbound}
y^a V_{1/2}^{(a)} (y) &= \delta_{0,a} + O_{\eps, \alpha} ( y^\alpha)  \\
y^a V_{1/2}^{(a)} (y) &\ll_A y^{-A}
}
for any integer $a \ge 0$, where $\delta_{0,a}$ is the Kronecker delta function. 
\end{rem}
\begin{proof}
By the residue theorem, 
\eqs{ 
\Lambda (s, \chi) = \frac 1 {2 \pi i} \int_R G(u) \Lambda(s+u,\chi) Y^u  \frac{du}{u}, }
where we use counterclockwise orientation over the rectangle $R$ with vertices $\pm 2\pm iT$ and $T>0$ does not correspond to the ordinate of a zero of $G(u)\Lambda(s+u, \chi)$. By Stirling's formula for $\Gamma(s)$ and the classical convexity bound, it follows that the horizontal integrals vanish as $T \to \infty$. For the vertical integral at $\re(u) = 2$, we have

\eqs{ 
\frac 1 {2\pi i} \int_{(2)} G(u) (2 \pi )^{-s-u} (3\cond{\chi})^{(s+u)/2} \Gamma(s+u) \biggl( \sum_{\fb} \frac{\chi(\fb)}{\N(\fb)^{s+u}}  \biggr) \frac{Y^u}{u} du  
} 
while for the vertical integral at $\re(u) = -2$ we change $u$ to $-u$ and use the functional equation to get  
\bs{ &\frac 1 {2\pi i} 
\int_{(2)} G(u) \Lambda(s-u, \chi) Y^{-u} \frac{du}  u  = \frac{W(\chi)}{\sqrt{\cond{\chi}}} \frac 1 {2\pi i} 
\int_{(2)} G(u) \Lambda(1-s+u, \overline{\chi}) Y^{-u} \frac{du}  u \\
& = \frac{W(\chi)}{\sqrt{\cond{\chi}}} \frac 1 {2\pi i} 
\int_{(2)} G(u) (2 \pi )^{s-u-1} (3\cond{\chi})^{(1-s+u)/2} \Gamma(1-s+u) \\
& \qquad \cdot \biggl( \sum_{\fb} \frac{\overline{\chi(\fb)}}{\N(\fb)^{1-s+u}}  \biggr) Y^{-u} \frac{du}  u.
}
The result follows by combining the two integrals and dividing both sides by $(3\cond{\chi})^{s/2} (2\pi)^{-s} \Gamma(s)$. 
\end{proof}

\subsection{First Mollified Moment}
Using the definition of the mollifier in \eqref{mollifier} we see that
\eqs{
M(c, \kappa) = \prod_{j=0}^J \su{\fa \in \cI_j\\ 0 \le \Omega(\fa) \le \ell_j} \frac{\chi_c (\fa)\lambda(\fa)f(\fa, J) \nu(\fa) }{\N(\fa)^{1/2}\kappa^{\Omega(\fa)}}  = \su{\fa} \frac{\chi_c (\fa)\lambda(\fa) f(\fa, J)\upsilon_J (\fa) } { \N(\fa)^{1/2}\kappa^{\Omega(\fa)}}, }
where
\eqs{
\upsilon_J (\fa) =  \su{\fa_j \in \cI_j \text{ for } 0 \le j \le J\\
\fa = \fa_0 \cdots \fa_J \\ 
0 \le  \Omega(\fa_j) \le \ell_j} \nu(\fa_0) \cdots \nu(\fa_J).}
The first mollified moment we shall study is given by
\eqn{\label{1stMoment}
\sum_{c \in \sF(X)} L(1/2, \chi_c) M(c,1) = \su{\fa} \frac{\lambda(\fa) f(\fa, J) \upsilon_J (\fa)}{\N(\fa)^{1/2}} \sum_{c \in \sF(X)} \chi_c (\fa) L(1/2, \chi_c). 
}

In Lemma \ref{lem:approxfnceqn} taking $s=1/2, \chi = \chi_c, G(u) \equiv 1$, replacing $Y$ by $Y(3\cond{\chi})^{-1/2}$ and factoring each ideal $\fb$ as $\fb = (1-\omega)^{r_2} b\ring$ with $b \equiv 1 \mod 3$ and $r_2 \ge 0$  gives 
\bsc{\label{AppFncEqnatOnehalf}
L(1/2,\chi_c) &= \sum_{r_2 \ge 0} 3^{-r_2/2} \sum_{b \equiv 1 \mod 3}  \frac{\chi_c (b)}{\N(b)^{1/2}}  V \biggl( \frac {3^{r_2}\N(b)} Y \biggr)  \\
&\quad +  \frac{W(\chi_c)}{\sqrt{\cond{c}}}
\sum_{r_2 \ge 0} 3^{-r_2/2} \sum_{b \equiv 1 \mod 3} \frac{\overline{\chi_c (b)}}{\N(b)^{1/2}} V \biggl(  \frac {3^{r_2}\N(b)Y} {3\cond{c}}  \biggl), }
where 
\eqs{
V (y) = V_{1/2}(y) = \frac 1 {2\pi i} \int_{2-i\infty}^{2+i\infty} (2 \pi y)^{-u} \frac{\Gamma(1/2+u)}{\Gamma(1/2)} \frac{du}{u}.
}
Here we used the fact that $\chi_c (1-\omega) = 1$ whenever $c \equiv 1 \mod 9$ for $c\in\sF(X)$, which follows from \cite{IR}[Chapter 9, Theorem 1], Supplement to the Cubic Reciprocity Law.

First write $\fa = a\ring$ with $a \equiv 1 \mod 3$. Hence, $\chi_c (\fa\fb) = \chi_c (ab)$. By the cubic reciprocity law, $\chi_c (\fa\fb) = \chi_c (ab) = \chi_{ab} (c)$ unless $ab=1$. Next, replace $a, b$ by $na, nb$ where $n, a, b \equiv 1 \mod 3$ and $(a,b)=1$. Then, substituting \eqref{AppFncEqnatOnehalf} into \eqref{1stMoment} 
we have that
\eqs{
\eqref{1stMoment} = \su{\fa} \frac{\lambda(\fa) f(\fa, J) \upsilon_J (\fa)}{\N(\fa)^{1/2}} \bigl(
S_1 + S_2 \bigr) +  \sum_{a, n \equiv 1 \mod 3}  \frac{\lambda(na) f(na, J) \upsilon_J (na)}{\N(a)^{1/2} \N(n)} S_3 (n) 
}
where
\bsc{\label{S123}
S_1 &= \sum_{r_2 \ge 0} 3^{-r_2/2} \su{b \equiv 1 \mod 3\\ ab = \cube} \frac 1 {\N(b)^{1/2}} V \biggl( \frac {3^{r_2}\N(b)} Y \biggl) \su{c \in \sF(X)\\ (c,ab)=1} 1    \\
S_2 &= \sum_{r_2 \ge 0} 3^{-r_2/2} \su{b \equiv 1 \mod 3\\ ab \neq \cube} \frac 1 {\N(b)^{1/2}} V \biggl( \frac {3^{r_2}\N(b)} Y \biggl) \su{c \in \sF(X)\\ (c,ab)=1} \chi_{ab} (c)  \\
S_3 (n) &= \sum_{r_2 \ge 0} 3^{-r_2/2} \su{b \equiv 1 \mod 3\\ (b,a) = 1} \frac 1 {\N(b)^{1/2}} \sum_{c \in \sF(X)} \chi_c (ab^2n^3) \frac{W(\chi_c)}{\sqrt{\cond{c}}}  V \biggl(  \frac {3^{r_2}\N(bn)Y} {3\cond{c}}  \biggl).
}

\subsection{The main term $S_1$} \label{sec:S1}
First recalling that $V(y) = 1 + O(y^{1/6-\eps})$ (see \eqref{Vbound}) and trivially estimating the sum over $c$ by $X\log X$ we see that
\eqs{
S_1 = \sum_{r_2 \ge 0} 3^{-r_2/2} \su{b \equiv 1 \mod 3\\ ab = \cube} \frac 1 {\N(b)^{1/2}} \su{c \in \sF(X)\\ (c,ab)=1} 1   +  O \biggl(\frac{X^{1+\eps}}{Y^{1/6}} \su{b \equiv 1 \mod 3\\ ab = \cube} \frac 1 {\N(b)^{1/3+\eps}} \biggr)
}

When $ab = \cube$, using Lemma \ref{familysize} (with $n = ab$) gives 
\eqs{\su{c \in \sF(X)\\ (c,ab)=1} 1 = C_1 (ab) X \log X  + C_2 (ab) X + O_\eps  \Bigl(  	\bigl( X^{3/4+ \eps} (\N(ab))^\eps \Bigr).}
Using \eqref{Fpsi(s)} we have for any $\alpha \in \ring$ that 
\bs{ \frac{F_{\psi_0}' (1;\alpha)}{F_{\psi_0} (1;\alpha) } &= \su{\pi \equiv 1 \mod 3\\ \pi \nmid \alpha}  \frac{6\log \N(\pi) \bigl(1 -  \frac 1 {\N(\pi)}\bigr)}{\N(\pi)^2 \Bigl(1 - \frac 3{\N(\pi)^2} + \frac 2 {\N(\pi)^3}\Bigr)} + 2 \su{\pi \equiv 1 \mod 3\\ \pi \mid \alpha} \frac{\log \N(\pi)}{\N(\pi) - 1}. }
Since $F_{\psi_0} (1;\alpha) \le 1$, it follows that $F_{\psi_0}' (1;\alpha) \ll (\N(\alpha))^\eps$, and thus $C_2(ab) \ll \N(ab)^\eps$. Therefore, recalling the definition of $C_1(ab)$ in \eqref{C12} we can now write
\bsc{\label{S1Main+Error}
S_1 &= X \log X  \frac {4\pi^2 \sqrt 3} {2187(\sqrt 3 - 1)}
\su{b \equiv 1 \mod 3\\ ab = \cube} \frac{F_{\psi_0}(1;ab)    }{\N(b)^{1/2}} \\
&\qquad +  O \biggl( \frac{X^{1+\eps}}{Y^{1/6}} \su{b \equiv 1 \mod 3\\ ab = \cube} \frac 1 {\N(b)^{1/3+\eps}} + X \su{b \equiv 1 \mod 3\\ ab = \cube} \frac {(\N(ab))^\eps}{\N(b)^{1/2}} \biggr).
}
Since $ab$ is a cube, writing $a = a_1 a_2^2 a_3^3$, where $a_1, a_2$ are square-free and co-prime, we can replace $b$ above by $a_1^2 a_2 b^3$ (with a different $b$) to get that the error term above is 
\eqn{\label{S1Error}
\ll 
\frac{X^{1+\eps}}{Y^{1/6}\N(a_1^2 a_2)^{1/3+\eps}}  + X \frac {\N(a_1a_2a_3)^\eps}{\N(a_1^2 a_2)^{1/2}}.
}

\subsection{Estimate of $S_2$}
Write $S_2$ as
\eqn{\label{S2integral}
S_2 = - \sum_{r_2 \ge 0} 3^{-r_2/2} \int_{1^-}^\infty  \frac {3^{r_2}} Y  V' \biggl( \frac {3^{r_2}z} Y \biggl) \su{b \equiv 1 \mod 3\\ ab \neq \cube\\ \N(b) \le z} \frac 1 {\N(b)^{1/2}}  \su{c \in \sF(X)\\ (c,ab)=1} \chi_{ab} (c) dz.}
We consider 
\bs{
\su{c \in \sF\\ \cond c \sim y} \chi_{ab} (c) &= \frac 1 {h_{(9)}} \sum_{\psi \mod 9} \su{d \equiv 1 \mod 3\\(d, ab)=1\\ \N(d) \le y^{1/2}
} \mu_K(d)  \psi (d^3) \\
& \quad \cdot
\ps{c_1 \equiv 1 \mod 3\\  (c_1, abd)=1} (\chi_{ab} \psi) (c_1)
\ps{c_2 \equiv 1 \mod 3\\ (c_2, abd)=1\\ \N(c_1c_2) \sim W} (\chi_{ab} \psi) (c_2^2)
}
for $y=X2^{-k} \ge 1$ with $k \ge 1$, where $W = y/\N(d)^2$. Here and in what follows, $\Sigma^\ast$ indicates summation over square-free integers. First we need to bound
\eqn{\label{S2dyadic}
\su{b \equiv 1 \mod 3\\ (ab,d) = 1 \\ ab \neq \cube\\ \N(b) \sim z}   \frac  1  {\N(b)^{1/2}} \ps{c_1 \equiv 1 \mod 3\\ (c_1, abd)=1 \\ \N(c_1) \sim U} (\chi_{ab} \psi) (c_1)
\ps{c_2 \equiv 1 \mod 3\\ \N(c_2) \sim V \\ (c_2, abd)=1 \\ \N(c_1c_2) \sim W} (\chi_{ab} \psi) (c_2^2)
}
for a fixed $d$ with $\N(d) \in [1,\sqrt y]$ and $1 \le U = W2^{-i}, V=W2^{-j}$ with $i, j \ge 1$ satisfying $W/4 < UV < 2W$. Applying Perron's formula with $T = W$ for the sums over $c_1$ and $c_2$ and then summing over $b$, we see that \eqref{S2dyadic} is bounded by 
\eqn{
\label{S2AfterPerron}
W \log W \sup_t	\su{b \equiv 1 \mod 3\\ (ab,d) = 1 \\ ab \neq \cube \\  \N(b) \sim z}   \frac 1 {\N(b)^{1/2}}  \big|\Sigma_1 (U) \Sigma_2 (V) \big|
+ O\bigl( W^\eps z^{1/2} \bigr),
}
where 
\eqs{
\Sigma_1 (U) = \ps{c_1 \equiv 1 \mod 3\\ (c_1, abd)=1 \\ \N(c_1) \sim U} \frac{(\chi_{ab} \psi) (c_1)}{\N(c_1)^s}, \qquad
\Sigma_2 (V) = \ps{c_2 \equiv 1 \mod 3\\(c_2, abd)=1 \\ \N(c_2) \sim V} \frac{(\chi_{ab} \psi) (c_2^2)}{\N(c_2)^s}
}
and $s= 1 + 1/\log 2W + it$. Without loss of generality, we can assume that $U \le V$ so that $V \gg W^{1/2}$. Removing the square-free condition on $c_2$, we have that
\bs{\Sigma_2 (V) &= \su{e \equiv 1 \mod 3\\(e, abd)=1 \\ \N(e) \le V^{1/2} } \frac{\mu_K(e) (\chi_{ab} \psi) (e^4)}{\N(e)^{2s}}  \su{c_2 \equiv 1 \mod 3\\(c_2, abd)=1 \\ \N(c_2e^2) \sim V} \frac{(\chi_{ab} \psi) (c_2^2)}{\N(c_2)^s}.
}
For the inner sum we use the estimate 
\eqs{\su{c_2 \equiv 1 \mod 3\\(c_2, abd)=1 \\ \N(c_2e^2) \sim V} \frac{(\chi_{ab} \psi) (c_2^2)}{\N(c_2)^s} \ll (V/\N(e)^2)^{-1/2+\eps} \N(ab)^{1/4+\eps},
}
which follows easily using Perron's formula and classical convexity bound together with partial integration, so that
\eqs{
\Sigma_2 (V) \ll V^{-1/2+\eps}\N(ab)^{1/4} \N(abd)^\eps. 
}
Using this result in \eqref{S2AfterPerron} and applying cubic large sieve inequality \cite[Theorem 2]{HB2000} for $\Sigma_1(U)$ yields that
\eqs{
\eqref{S2AfterPerron} \ll \N(a)^{1/4} \N(ad)^\eps \bigl(W^{ 1/2+\eps} z^{3/4+3\eps/2}  +  W^{3/4+\eps} z^{1/4+2\eps} +  W^{2/3+\eps}  z^{7/12+\eps}\bigr).  
}
Summing this over all $\ll \log W$ possible $U$'s and then over $\N(d) \le y^{1/2}$, and finally over $\ll \log X$ possible $y$'s and over $\ll \log z$ possible $z$'s, and inserting the resulting estimates in \eqref{S2integral} we conclude that 
\eqn{\label{S2bound}
S_2 \ll \N(a)^{1/4+\eps} \bigl(X^{ 1/2+\eps} Y^{3/4+\eps}  +  X^{3/4+\eps} Y^{1/4+\eps} + X^{2/3+\eps}  Y^{7/12+\eps}\bigr),
}
upon redefining $\eps$.
\subsection{Estimate of $S_3(n)$}
Since $V(y) \ll y^{-1}$ we see that the contribution of $c$ with $\cond c \le Y$ to $S_3(n)$ is $\ll Y\log Y/\N(n)$. Hence it is enough to consider 
\eqs{
\su{b \equiv 1 \mod 3\\ (b,a) = 1} \frac 1 {\N(b)^{1/2}} \su{c \in \sF\\ Y < \cond c \le X} \chi_c (ab^2n^3) \frac{W(\chi_c)}{\sqrt{\cond{c}}}  V \biggl(  \frac {3^{r_2}\N(bn)Y} {3\cond{c}}  \biggl).}
Using partial integration twice expresses this sum as
\bsc{\label{S3integral}
&= - \int_{1^-}^\infty  \frac Z X V' \biggl( \frac {Zz} X \biggr) \su{b \equiv 1 \mod 3\\ (b,a)=1\\ \N(b) \le z}  \frac 1 {\N(b)^{1/2}} \su{c \in \sF\\ Y < \cond c \le X} \chi_c (ab^2n^3) \frac {W(\chi_c)}{\sqrt{\cond{c}}} dz \\
& \quad - \int_{1^-}^\infty \int_Y^X  \biggl( \frac Z {y^2} V' \biggl( \frac {Zz} y \biggr) + \frac {Z^2z} {y^3} V'' \biggl( \frac {Zz} y \biggr)\biggr)  \\
& \qquad \cdot 
\su{b \equiv 1 \mod 3\\ (b,a)=1\\ \N(b) \le z}  \frac 1 {\N(b)^{1/2}}  \biggl(\su{c \in \sF\\ Y < \cond c \le y} \chi_c (ab^2n^3) \frac {W(\chi_c)}{\sqrt{\cond{c}}} \biggr)  dy dz,
}
where $Z = 3^{r_2-1} Y\N(n)$ . First note that for $c = c_1^2c_2 \in \sF$, 
\eqs{W(\chi_c) = \chi_{c_1^2} (c_2) \chi_{c_2} (c_1) \overline{W(\chi_{c_1})} W(\chi_{c_2}) = \overline{W(\chi_{c_1})} W(\chi_{c_2}),}
where the second equality follows by cubic reciprocity law. Note also that 
\eqs{W(\chi_{c_i}) = \chi_{c_i} (\sqrt{-3}) g(1,c_i).}
Hence, 
\eqs{W(\chi_c) =  \chi_c (\sqrt{-3}) \overline{g(1,c_1)} g(1,c_2) = \overline{g(1,c_1)} g(1,c_2),}
since $\chi_c (\sqrt{-3}) = \chi_{c} (w(1-\omega)) = 1$ whenever $c \equiv 1 \mod 9$. It follows that 
\bs{
&\su{c \in \sF\\ \cond c \sim y} \chi_c (ab^2n^3) \frac{W(\chi_c)}{(\cond{c})^{1/2}} 
\\
&= \frac 1 {h_{(9)}} \sum_{\psi \mod 9} \su{c_1 \equiv 1 \mod 3\\  (c_1, abn)=1}  \psi(c_1^2)  \frac{\overline{g(a^2bn^3,c_1)}}{\N(c_1)^{1/2}} 
\su{c_2 \equiv 1 \mod 3\\ (c_2, abnc_1)=1\\ \N(c_1c_2) \sim y} \psi(c_2)	\frac{g(a^2bn^3,c_2)}{\N(c_2)^{1/2}},
}
where 
we used  Lemma \ref{cubic-GS-lemma1} and \ref{cubic-GS-lemma2}  to replace $\chi_{c_2} (ab^2n^3)g(1,c_2)$ by $g(a^2bn^3, c_2)$, $\chi_{c_1^2} (ab^2n^3) \overline{g(1,c_1)}$ by $\overline{g(a^2bn^3,c_1)}$, and also to remove the condition that $c_1, c_2$ be square-free. Next we remove the condition that $(c_1, c_2)=1$ to get 
\bs{
&= \frac 1 {h_{(9)}} \sum_{\psi \mod 9} \su{d \equiv 1 \mod 3\\(d, abn) = 1 \\ \N(d)^2 \le 2y} \mu_K(d)\psi(d^3) \\
&\quad \cdot \su{c_1 \equiv 1 \mod 3\\  (c_1, abdn)=1}  \psi(c_1^2)  \frac{\overline{g(a^2bdn^3,c_1)}}{\N(c_1)^{1/2}} 
\su{c_2 \equiv 1 \mod 3\\ (c_2, abdn)=1\\ \N(c_1c_2) \sim y/\N(d)^2} \psi(c_2)	\frac{g(a^2bdn^3, c_2)}{\N(c_2)^{1/2}}
}
where we used Lemma \ref{cubic-GS-lemma1} to first introduce the conditions $(d,c_1c_2)=1$, which can be done since otherwise the Gauss sums vanish, and then write $g(a^2bn^3,dc_i)$ as $g(a^2bn^3d, c_i)g(a^2bn^3, d)$, and finally note that $|g(a^2bn^3,d)|^2 = \N(d)$ as $(d,abn)=1$. We shall first estimate
\eqn{\label{S3dyadicaverage}
\su{b \equiv 1 \mod 3\\ (b,ad)=1\\ \N(b) \sim z}  \frac 1 {\N(b)^{1/2}} \su{c_1 \equiv 1 \mod 3\\  (c_1, abdn)=1\\ \N(c_1) \sim U}  \psi(c_1^2)  \frac{\overline{g(a^2bdn^3,c_1)}}{\N(c_1)^{1/2}} \su{c_2 \equiv 1 \mod 3\\ (c_2, abdn)=1\\ \N(c_2) \sim V \\ \N(c_1c_2) \sim W} \psi(c_2)	\frac{g(a^2bdn^3, c_2)}{\N(c_2)^{1/2}} 
}
for a fixed $d$ with $\N(d) \in [1,\sqrt y]$ and $1 \le U = W2^{-i}, V= W2^{-j}$ satisfying $W/4 < UV < 2W$, where $W= y/\N(d)^2$. Applying Perron's formula with $T = W$, we see that \eqref{S3dyadicaverage} is bounded by 
\eqn{
\label{LargeD}
W \log W \sup_t	\su{b \equiv 1 \mod 3\\ (b,ad) = 1 \\ \N(b) \sim z}   \frac 1 {\N(b)^{1/2}}  \big|\Sigma_1 (U) \Sigma_2 (V) \big|
+ O\bigl( W^\eps z^{1/2} \bigr)
}
where 
\eqs{
\Sigma_1 (U) = \su{c_1 \equiv 1 \mod 3\\  (c_1, abdn)=1\\ \N(c_1) \sim U}  \psi(c_1^2)  \frac{\overline{g(a^2bdn^3,c_1)}}{\N(c_1)^{1/2+s}}, \quad \Sigma_2 (V) = \su{c_2 \equiv 1 \mod 3\\ (c_2, abdn)=1\\ \N(c_2) \sim V} \psi(c_2)	\frac{g(a^2bdn^3, c_2)}{\N(c_2)^{1/2+s}} 
}
and $s= 1 + 1/\log 2W + it$. By the last equation in the proof of \cite[Proposition 6.2]{DG}
\multn{\label{S3PV}
\su{c \equiv 1 \mod 3\\  (c, r)=1\\ \N(c) \le z} \lambda(c) \frac{g(r, c)}{\N(c)^{1/2}} \\
\ll_{\cond{\lambda}} \N(r_1r_3^\ast)^\eps \bigl( z^{5/6}\N(r_1)^{-1/6} + z^{2/3+\eps} \N(r_1r_2^2)^{1/6} + z^{1/2+\eps} \N(r_1r_2^2)^{1/4} \bigr),	}
where $r=r_1r_2^2r_3^3$ with $r_i \equiv 1 \mod 3$, $r_1, r_2$ coprime and square-free and $r_3^\ast$ is the product of primes dividing $r_3$ but not $r_1r_2$. Also, the first error term appears only when $r_2 = 1$.

Now, writing $a= a_1a_2^2a_3^3$ and $b= b_1b_2^2b_3^3$ we see that $a^2bdn^3 = r_1 r_2^2 r_3^3$ with $r_1 = a_2b_1d, r_2 = a_1 b_2$ and $r_3 =b_3a_2 n$. Since $(abn,d)=1$ and $d$ is square-free, it follows using partial integration and \eqref{S3PV}  that 
\eqs{
\Sigma_2 (V) \ll \N(abdn)^\eps \Bigl( 
\frac{1}{V^{1/6}\N(a_2b_1d)^{1/6}} +  \frac{\N(a_1^2a_2 d b_1b_2^2 )^{1/6}}{V^{ 1/3-\eps}} + \frac{\N(a_1^2a_2 d b_1b_2^2 )^{1/4}}{V^{ 1/2-\eps}}
\Bigr)  
}
where the first term appears only when $a_1  = b_2 = 1$.  Using this estimate we see that the sum over $b$ in \eqref{LargeD} is bounded by 
\bs{
& \frac{\N(adn)^\eps}{z^{1/6-\eps}V^{1/6} \N(a_2d)^{1/6}} \sum_{\N(b_3) \le z^{1/3}} \N(b_3)^{-1} \biggl( \ps{b_1 \equiv 1 \mod 3\\ \N(b_1) \sim z/\N(b_3)^3} |\Sigma_1(U)|^2\biggr)^{1/2} \\
& + \frac{\N(a_1^2 a_2 d)^{1/6} \N(adn)^\eps z^{1/6+\eps}}{V^{ 1/3-\eps}}  \sum_{\N(b_3) \le z^{1/3}} \N(b_3)^{-2} \\
&\qquad \cdot \ps{\N(b_2^2b_3^3) \le z} \N(b_2)^{-1} \biggl( \; \ps{b_1 \equiv 1 \mod 3\\ \N(b_1) \sim z/\N(b_2^2b_3^3)} |\Sigma_1(U)|^2\biggr)^{1/2} \\
& + \frac{\N(a_1^2 a_2 d)^{1/4} \N(adn)^\eps z^{1/4+\eps} }{V^{ 1/2-\eps}}  \sum_{\N(b_3) \le z^{1/3}} \N(b_3)^{-9/4} \\
& \qquad \cdot \ps{\N(b_2^2b_3^3) \le z} \N(b_2)^{-1} \biggl( \; \ps{b_1 \equiv 1 \mod 3\\ \N(b_1) \sim z/\N(b_2^2b_3^3)} |\Sigma_1(U)|^2\biggr)^{1/2}, 
}
where the first term appears only when $a_1 = 1$. Now, applying cubic large-sieve inequality for the innermost sums over $b_1$ and then summing over $b_2$ and $b_3$ gives the bound
\eqs{
U^{-1/2+\eps/2} z^{1/2+3\eps/2} + U^{\eps/2} z^{3\eps/2} + z^{1/3+3\eps/2} U^{-1/6+\eps/2}
}
for each of the three terms. Now using this estimate and assuming that $U \le V$, we conclude that
\bs{
\eqref{LargeD} 
&\ll \N(adn)^\eps \Bigl(
\frac{W^{5/6+\eps} }{\N(a_2d)^{1/6}} z^{1/3+5\eps/2} + \frac{W^{11/12+\eps} }{\N(a_2d)^{1/6}} z^{5\eps/2-1/6} \\
& + \frac{W^{5/6+\eps} }{\N(a_2d)^{1/6}} z^{1/6+5\eps/2} + W^{2/3+\eps} \N(a_1^2 a_2 d)^{1/6}   z^{2/3+5\eps/2} \\
& + W^{5/6+\eps}\N(a_1^2 a_2 d)^{1/6}  z^{1/6+5\eps/2} + W^{3/4+\eps}\N(a_1^2 a_2 d)^{1/6}  z^{1/2+5\eps/2} 	\\
& + \quad W^{1/2+\eps} \N(a_1^2 a_2 d)^{1/4}  z^{3/4+5\eps/2}  + W^{3/4+\eps} \N(a_1^2 a_2 d)^{1/4}  z^{1/4+5\eps/2} \\
& + W^{2/3+\eps} \N(a_1^2 a_2 d)^{1/4}  z^{7/12+5\eps/2} \Bigr) + W^\eps z^{1/2}, }
where the first three terms appear only when $a_1 = 1$. Summing up over the range $\N(d) \le B$ yields the bound
\bsc{	\label{S3smalld} 
\N(an)^\eps \Bigl( &
\frac{y^{5/6+\eps} }{\N(a_2)^{1/6}} z^{1/3+5\eps/2} + \frac{y^{11/12+\eps} }{\N(a_2)^{1/6}} z^{5\eps/2-1/6} \\
& + y^{2/3+\eps} \N(a_1^2 a_2 )^{1/6}   z^{2/3+5\eps/2} + y^{5/6+\eps}\N(a_1^2 a_2 )^{1/6}  z^{1/6+5\eps/2} \\
& + y^{3/4+\eps}\N(a_1^2 a_2 )^{1/6}  z^{1/2+5\eps/2} + y^{3/4+\eps} \N(a_1^2 a_2 )^{1/4}  z^{1/4+5\eps/2} 	\\
& + y^{1/2+\eps} \N(a_1^2 a_2 )^{1/4}  z^{3/4+5\eps/2} B^{1/4} + y^{2/3+\eps} \N(a_1^2 a_2)^{1/4}  z^{7/12+5\eps/2} \Bigr), }
where  $a_1 = 1$ for the first two terms.

Assume now that $B < \N(d) \le (2y)^{1/2}$. Write $b=b_1b_2^2$ with $b_i \equiv 1 \mod 3$ and $b_1$ square-free. Then, Cauchy-Schwarz inequality yields that
\eqs{
\su{b \equiv 1 \mod 3\\ (b,ad) = 1 \\ \N(b) \sim z} \frac {\big|\Sigma_1 (U) \Sigma_2 (V) \big|} {\N(b)^{1/2}}   
\le z^{-1/2}\su{b_2 \equiv 1 \mod 3\\ \N(b_2) \le \sqrt z} 	\biggl( \; \ps{b_1 \equiv 1 \mod 3\\ \N(b_1b_2^2) \sim z}    \big|\Sigma_1 (U) \big|^2 \ps{b_1 \equiv 1 \mod 3\\ \N(b_1b_2^2) \sim z}   \big| \Sigma_2 (V) \big|^2 \biggr)^{1/2}.
}
Applying cubic large sieve inequality for both sums we conclude that
\mult{
\su{b \equiv 1 \mod 3\\ (b,ad) = 1 \\ \N(b) \sim z}   \frac {\big|\Sigma_1 (U) \Sigma_2 (V) \big| } {\N(b)^{1/2}}   \ll W^{\eps/2} \biggl(
W^{-1/2} z^{1/2+\eps}   +  U^{-1/2} z^{\eps}   \\
+   U^{-1/2} z^{1/3+\eps} V^{-1/6}  +   V^{-1/2} z^{\eps}  +   1  +   z^{1/3+ \eps} U^{-1/6}  V^{-1/2} +   W^{-1/6} z^{1/6+\eps}  
\biggr).
}
Inserting this in \eqref{LargeD} and
assuming $V \le U$ shows that \eqref{S3dyadicaverage} is bounded by
\eqs{
W^{1/2+2\eps} z^{1/2+\eps}   +  W^{3/4+2\eps} z^{\eps}  +  W^{2/3+2\eps} z^{1/3+\eps} +  W^{1+\eps} z^{\eps}  +   z^{1/3+ \eps} W^{5/6+\eps}.
}
Summing over $B < \N(d) \le \sqrt {2y}$ we get the bound	
\multn{\label{S3Larged}
y^{1/2+2\eps} z^{1/2+\eps}   +  y^{3/4+2\eps} z^{\eps} B^{-1/2} +  y^{2/3+2\eps} z^{1/3+\eps} B^{-1/3} \\
+  y^{1+\eps} z^{\eps} B^{-1} +    y^{5/6+\eps}  z^{1/3+ \eps} B^{-2/3}.
}

Combining \eqref{S3smalld} and \eqref{S3Larged}, balancing terms using Lemma \ref{balancing} with $B \in [1,(2y)^{1/2}]$ shows that for $y \in (Y,X]$, 
\eqs{ \su{b \equiv 1 \mod 3\\ (b,a)=1\\ \N(b) \le z}  \frac 1 {\N(b)^{1/2}}  \biggl(\su{c \in \sF\\ Y < \cond c \le y} \chi_c (ab^2n^3) \frac {W(\chi_c)}{\sqrt{\cond{c}}} \biggr)  
}
is bounded by
\bsc{\label{S3rangescombined}
(yz&\N(an))^\eps \Bigl( \frac{y^{5/6}z^{1/3} }{\N(a_2)^{1/6}}  + \frac{y^{11/12} }{\N(a_2)^{1/6} z^{1/6}}  + y^{2/3} \N(a_1^2 a_2 )^{1/6}   z^{2/3} \\
&	+ y^{2/3} \N(a_1^2 a_2)^{1/4}  z^{7/12} + y^{5/6}\N(a_1^2 a_2 )^{1/6}  z^{1/6}  
+ y^{3/4}\N(a_1^2 a_2 )^{1/6}  z^{1/2} \\
&+ y^{1/2} \N(a_1^2 a_2 )^{1/4}  z^{3/4}
+ y^{3/4} \N(a_1^2 a_2 )^{1/4}  z^{1/4} 
+  y^{3/5} \N(a_1^2 a_2)^{1/5}  z^{3/5}\\
&  + y^\eps z^{3\eps}	y^{13/22}  z^{7/11} \N(a_1^2 a_2 )^{2/11} \Bigr) ,  }
upon redefining $\eps$, where  $a_1 = 1$ for the first two terms. 

Note that
\eqs{
\bigg|\sum_{a, n \equiv 1 \mod 3}  \frac{\lambda(na) f(na, J) \upsilon_J (na)}{\N(a)^{1/2} \N(n)} S_3 (n) \bigg| \le \sum_{a\equiv 1 \mod 3}  \frac{1}{\N(a)^{1/2}} S_3, 
}
where
\eqs{S_3 = \sum_{n \equiv 1 \mod 3}  \frac {\upsilon_J (an)} {\N(n)} |S_3 (n)|.}
We split the sum over $n$ at $X/Y$, and for $\N(n) \le X/Y$, we also split the $z$-integral in \eqref{S3integral} at $X/(Y\N(n))$. Using the bounds in \eqref{Vbound}, we see that the contribution of each term in \eqref{S3rangescombined} of the form $y^u z^v$ with $u,v \in (0,1)$ to $S_3(n)$ will be 
\eqs{X^{u+v} Y^{-v} \sum_n \frac{\upsilon_J (an)}{\N(n)^{1+v -\eps}}}
The second term in \eqref{S3rangescombined} contributes 
\eqn{\label{S31}
S_{31} := \frac{X^{11/12+\eps} }{\N(a_2)^{1/6}} 
\sum_n \frac{\upsilon_J (an) \N(an)^\eps }{\N(n)^{1 -\eps}},
}
where $a_1 = 1$. Observe that since $\upsilon_J (na) \le \upsilon_J (n) \upsilon_J (a)$,
\eqn{\label{convergentupsilonsum}
\sum_n \frac{\upsilon_J (an)}{\N(n)^{1+v -\eps}} < \upsilon_J (a) \exp \Bigl( \sum_{k_0 < \N\fp \le X^\theta_J} \frac 1 {\N(\fp)^{1+v-\eps}} \Bigr) \ll \upsilon_J (a).
}
Therefore, upon redefining $\eps$, we obtain
\bs{
S_3	&\ll S_{31}  +  \upsilon_J (a) S_{32},
}
where
\bsc{\label{S32}
(X\N(a))^{-\eps} S_{32} &= 
\frac{X^{7/6}}{\N(a_2)^{1/6} Y^{1/3}} + \frac{X^{4/3}  \N(a_1^2 a_2 )^{1/6} }{Y^{ 2/3}}+ \frac{X^{5/4}  \N(a_1^2 a_2)^{1/4} }{Y^{ 7/12}} \\
& +  \frac{X\N(a_1^2 a_2 )^{1/6}}{Y^{1/6}} + \frac{X^{5/4} \N(a_1^2 a_2 )^{1/6}}{Y^{ 1/2}}	
+ \frac{X^{5/4}\N(a_1^2 a_2 )^{1/4}}{Y^{ 3/4}} \\
& + \frac{X  \N(a_1^2 a_2 )^{1/4} }{Y^{ 1/4}} 
+ \frac{X^{6/5}\N(a_1^2 a_2)^{1/5}  }{Y^{ 3/5} }
+ \frac{X^{27/22}  \N(a_1^2 a_2 )^{2/11}}{ Y^{ 7/11}}
}
and $a_1 = 1$ for the first term. 

\subsection{Gluing the pieces together}
Combining \eqref{S1Main+Error}, \eqref{S1Error}, \eqref{S2bound} and \eqref{S32} and inserting in \eqref{1stMoment} we get
\bs{
\eqref{1stMoment} 
&- \frac {4\pi^2 \sqrt 3 X \log X } {2187(\sqrt 3 - 1)} \sum_{a} \frac{\lambda(a) f(a, J) \upsilon_J (a)}{\N(a)^{1/2}}  
\su{b \equiv 1 \mod 3\\ ab = \cube} \frac{F_{\psi_0}(1;ab)    }{\N(b)^{1/2}} \\
& \ll \Sigma_1 + \Sigma_2 + \Sigma_3,
}
where
\eqs{\Sigma_1 = X \ps{a_1,a_2 \equiv 1 \mod 3\\ (a_1,a_2)=1} \sum_{a_3 \equiv 1 \mod 3} \frac{\upsilon_J (a_1a_2^2a_3^3)}{\N(a_1a_2a_3)^{3/2-\eps}},}
comes from the second error term in \eqref{S1Error}, 
\bs{\Sigma_2 &= X^{11/12+\eps} \sum_{a\equiv 1 \mod 3}  \frac{1}{\N(a)^{1/2}} 
\sum_n \frac{\upsilon_J (an) \N(an)^\eps }{\N(n)^{1 -\eps}} \\
&\ll X^{11/12+\eps} \sum_b 
\upsilon_J (b) \N(b)^{2\eps-1/2}
}
results from $S_{31}$, and
\eqs{\Sigma_3 = \sum_{a\equiv 1 \mod 3}  \frac{\upsilon_J (a)}{\N(a)^{1/2}}  (S_{32} + |S_2|).}

Proceeding as in \eqref{convergentupsilonsum} we can see that $\Sigma_1 \ll X$.

Note that the first error term in \eqref{S1Error} is smaller than the fourth term of $S_{32}$ in \eqref{S32}, thus was omitted above. Next, combining the bounds \eqref{S2bound} and \eqref{S32} for $S_2$ and $S_{32}$ and balancing terms using Lemma \ref{balancing} with $Y \in [1,X]$ shows that 
\eqs{S_2 + S_{32} \ll (X\N(a))^\eps \sum_i X^{u_i} \N(a_1)^{t_{1i}} \N(a_2)^{t_{2i}} \N(a_3)^{t_{3i}} }
where $(u_i, t_{1i}, t_{2i}, t_{3i})$ are given by
\begingroup\makeatletter\def\f@size{8}\check@mathfonts
\bs{
&(3/4, 1/4, 1/2, 3/4),
(3/4, 1/2, 1/4, 0),
(7/8, 5/16, 7/16, 9/16),\\
&(7/8, 7/16, 5/16, 3/16),
(15/17, 5/17, 7/17, 9/17),(9/10, 13/40, 17/40, 21/40),\\
&(9/10, 17/40, 13/40, 9/40),
(10/11, 3/11, 9/22, 6/11),
(11/12, 5/18, 7/18, 1/2),\\
&(59/64, 23/64, 25/64, 27/64),
(59/64, 25/64, 23/64, 21/64),
(13/14, 1/7, 5/14, 3/7),\\
&(66/71, 23/71, 25/71, 27/71),
(43/46, 7/23, 8/23, 9/23), 
(23/24, 3/8, 3/8, 3/8),\\
&(44/45, 13/45, 31/90, 2/5),
(51/52, 23/78, 25/78, 9/26),
(65/66, 1/11, 19/66, 3/11).
}
\endgroup

Summing over $a$, we therefore conclude that  
\bs{\Sigma_3 &\ll  \sum_i X^{u_i+\eps} \ps{a_1 \equiv 1 \mod 3}  \upsilon_J (a_1) N(a_1)^{\eps+t_{1i}-1/2} \\
& \cdot \ps{a_2 \equiv 1 \mod 3}  \frac{\upsilon_J(a_2^2) }{ N(a_2)^{1-2\eps-t_{2i}}}
\sum_{a_3 \equiv 1 \mod 3}  \frac{\upsilon_J(a_3^3)}{N(a_3)^{3/2-3\eps-t_{3i}}}.   
}
Now, we claim that $\Sigma_3 \ll \sum_i X^{v_i}$, where
\eqn{\label{claimsigma3}
v_i = 
\begin{cases}
u_i + \eps + \eta(1/2 + t_{1i} + t_{2i}/2 + 2\eps) & \text{if } t_{3i} < 1/2, \\ 
u_i + \eps + \eta(1/3 + t_{1i} + t_{2i}/2 + t_{3i}/3 + 3\eps) & \text{otherwise,}	
\end{cases}
}
and
\eqs{
\eta = \sum_{0\le r\le J} \theta_r \ell_r \le \frac{2e^{1-\beta}  }{e^{1-\beta}-1} \Theta^{1-\beta}. 
}
Assuming this claim holds for the moment, we see that $\Sigma_3 \ll X$ if each $v_i \le 1$, which in turn holds provided that
\eqn{\label{ThetaCond1stmoment}
\Theta \le  \Bigl(\frac{e^{1-\beta}-1}{2e^{1-\beta}  } \min \Big\{ \frac{2 - 132 \eps}{264\eps + 97}, \frac{3-156\eps}{312\eps + 149} \Big \}  \Bigr)^{\frac 1 {1-\beta}}.
}

To prove the claim \eqref{claimsigma3}, we first note that for $t_{3i}<1/2$, it follows as in \eqref{convergentupsilonsum} that the sum over $a_3$ is $\ll 1$. For the rest of the sums over $a_1, a_2$ and $a_3$ we have
\bs{
\sum_{a_j} \frac{\nu(a_j^j)}{\N(a_j)^{\alpha_j} } 
&= \prod_{0 \le r \le J} \su{a_{rj} \in \cI_r\\ 0 \le \Omega(a_r) \le \ell_r/j} \frac{\nu(a_{rj}^j)}{\N(a_{rj})^{\alpha_j}} \\
& = \prod_{0 \le r \le J}  E_{\ell_r/j} \Bigl(\sum_{\fp \in \cI_r}  (\N\fp)^{-\alpha_j} \Bigr),
}
where $j=1,2,3$ and $\alpha_j$ represents the corresponding exponents. Note that 
\bs{
\sum_{\fp \in \cI_r}  (\N\fp)^{-\alpha_j} 
&< 2\int_1^{X^{\theta_r/2}} t^{-\alpha_j} dt + \int_{X^{\theta_r/2}}^{X^{\theta_r}} t^{-\alpha_j} d\pi (t) 
< X^{\theta_r(1-\alpha_j),}
}
provided that
\eqs{
\frac{\log X}{(\log_2 X)^\alpha} > \max \Big\{\log 9, \frac{12}{1- \alpha_j}, \frac 2 {1-\alpha_j}\log \frac 4{1-\alpha_j} \Big\}.
}
Here, we used the inequality
\eqs{\pi(t) < 3t/\log t, \qquad (t>3). } 
If we further assume that
\eqs{
\frac{\log X}{(\log_2 X)^\alpha \log_3 X} > \frac{\alpha\beta}{1-\alpha_j}
}
holds for all $X \ge X_0$ for some $X_0$, and
\eqs{\Theta < 5^{-1/\beta}, }
then it follows from Remark \eqref{Elbound} that for $\alpha_j < 1$ and $0 \le r \le J$, 
\eqs{
E_{\ell_r/j} \Bigl( X^{\theta_r(1-\alpha_j)} \Bigr) < X^{\theta_r\ell_r (1-\alpha_j)/j}.
}
With these assumptions, we see that for $\alpha_j < 1$,
\eqs{
\sum_{a_j} \frac{\nu(a_j^j)}{\N(a_j)^{\alpha_j} } < X^{\eta(1-\alpha_j)/j},
}
Combining these estimates for the three sums over $a_1, a_2$ and $a_3$ gives the claim. 

Finally note that $\Sigma_2 \ll X$ holds as long as \eqref{ThetaCond1stmoment} is satisfied and $X$ is sufficiently large in terms of $\eps$ and $\alpha$.

\subsection{Proof of Theorem \ref{thm:firstmoment}}
Rewrite $F_{\Psi_0}(1, ab)$ as 
\eqs{
\pr{\pi \equiv 1 \mod 3}  \Bigl( 1 - \frac 3 {\N(\pi)^2} + \frac 2 {\N(\pi)^3 }\Bigr) \pr{\pi \equiv 1 \mod 3\\ \pi \mid \fa b} \frac {\N(\pi) } {\N(\pi) + 2}. 
}
Then, the main term of the first mollified moment is given by 
\eqs{ 
c_0 X\log X \su{a \equiv 1 \mod 3} \frac{\lambda(a) f(a, J) \upsilon_J (a) }{\N(a)^{1/2}} 
\su{b \equiv 1 \mod 3\\ ab = \cube} \frac 1  {\N(b)^{1/2}} \pr{\pi \equiv 1 \mod 3\\ \pi \mid \fa b} \frac {\N(\pi) } {\N(\pi) + 2}, 
}
where
\eqn{\label{c0}
c_0 = \frac {4\pi^2 \sqrt 3 } {2187(\sqrt 3 - 1)} \pr{\pi \equiv 1 \mod 3}  \Bigl( 1 - \frac 3 {\N(\pi)^2} + \frac 2 {\N(\pi)^3 }\Bigr).}
Given $a = \prod_r a_r$ with $a_r \in \cI_r$ and $0 \le \Omega(a_r) \le \ell_r$, we see that $ab = \cube$ if and only if $b = b'^3 \prod_r b_r$ with $b_r \in \cI_r$ such that $a_r b_r = \cube$ for each $r$, and $b' \not\in \cI_r$ for any $r$. Therefore, the sum over $a$ above can be written as
\bs{
&\su{b \equiv 1 \mod 3\\ \fp \mid b \Rightarrow \N\fp > X^\theta_J} \frac 1  {\N(b)^{3/2}} \pr{\pi \equiv 1 \mod 3\\ \pi \mid b}  \frac {\N(\pi) } {\N(\pi) + 2} \\
&\cdot \prod_{0 \le r \le J} \su{a_r \in \cI_r\\ 0 \le \Omega(a_r) \le \ell_r } \frac{\lambda(a_r)f(a_r,J)\nu(a_r)}{\N(a_r)^{1/2}} \su{b_r \in \cI_r\\ a_rb_r = \cube} \frac 1 {\N(b_r)^{1/2}} \pr{\pi \equiv 1 \mod 3\\ \pi \mid a_r b_r}  \frac {\N(\pi) } {\N(\pi) + 2}.
}
The first sum over $b$ equals
\eqn{\label{maintermlowerbound1}
	\pr{\xi \equiv 1 \mod 3\\ \N\xi > X^{\theta_J}} \biggl(1 + \frac {\N(\xi)}{(\N(\xi) + 2)(\N(\xi)^{3/2}-1)} \biggr) > 1.
}
Furthermore, note that writing $a_r = a_{1r}a_{2r}^2 a_{3r}^3$ with $a_{ir} \equiv 1 \mod $, $a_{1r}, a_{2r}$ square-free and co-prime, we see that $a_r b_r = \cube$ if and only if $b_r = a_{2r} a_{1r}^2 (b'_r)^3$. Thus, the sum over $a_r$ for each $r$ is 
\bs{
& \ll \ps{a_{1r} \in \cI_r} \frac{ \nu(a_{1r})} {\N(a_{1r})^{3/2}} \ps{a_{21r} \in \cI_r} \frac{  \nu(a_{12}^2)} {\N(a_{2r})^{3/2}} \su{a_{3r} \in \cI_r} \frac{\nu(a_{3r}^3)} {\N(a_{3r})^{3/2}} \\
& = \prod_{\xi \in \cI_r} \Bigl( 1 + \frac 1 {\N(\xi)^{3/2}} \Bigr) \Bigl( 1 + \frac 1 {2\N(\xi)^{3/2}} \Bigr) \prod_{\xi \in \cI_r} \Bigl( 1 + \sum_{n \ge 1} \frac 1 {(3n)! \N(\xi)^{3n/2}}\Bigr) \ll 1.
}
Now, we use this observation to remove the condition that $\Omega(a_r) \le \ell_r$. Put $T = T_r = X^{\theta_{r-1}/12}$. Since $T>1$, the sum over $a_r$ with $\Omega(a_r) > \ell_r$ is then
\bs{
&< T^{-\ell_r} \su{a_{1r}, a_{2r}, a_{3r} \in \cI_r} \frac{ T^{\Omega(a_{1r}a_{2r}^2 a_{3r}^3)} \nu(a_{1r}a_{2r}^2 a_{3r}^3)} {\N(a_{1r}a_{2r}^2 a_{3r}^3)^{1/2}} \su{b_r \in \cI_r} \frac 1 {\N(a_{2r}a_{1r}^2b_r^3)^{1/2}} \\
& \ll T^{-\ell_r} \ps{a_{1r} \in \cI_r} \frac{ T^{\Omega(a_{1r})} \nu(a_{1r})} {\N(a_{1r})^{3/2}} \ps{a_{21r} \in \cI_r} \frac{ T^{2\Omega(a_{2r})} \nu(a_{12}^2)} {\N(a_{2r})^{3/2}} \su{a_{3r} \in \cI_r} \frac{ T^{3\Omega(a_{3r})} \nu(a_{3r}^3)} {\N(a_{3r})^{3/2}} \\
& = T^{-\ell_r} \prod_{\xi \in \cI_r} \Bigl( 1 + \frac{T}{\N(\xi)^{3/2}} \Bigr) \Bigl( 1 + \frac{T^2}{2\N(\xi)^{3/2}} \Bigr) \prod_{\xi \in \cI_r} \Bigl( 1 + \sum_{n \ge 1} \frac{T^{3n}}{(3n)! \N(\xi)^{3n/2}}\Bigr) \\
& \ll T^{-\ell_r},
}
since $T^{3} < \N(\xi)^{1/4}$ for $\xi \in \cI_r$. Therefore, removing the restriction on the prime factors of $a_r$ for one of the $J+1$ sums over $a_r$ introduces an error of size $\ll T_r^{-\ell_r}$ since the other $J$ sums are $\ll 1$, and so is the sum over $b$. Since $\ell_r > \theta_r^{-\beta}$ for sufficiently large $X$, we see that for $r \ge 1$,
\eqs{
X^{-\theta_{r-1}\ell_r/12} < X^{-\theta_{r}^{1-\beta}/12e} < X^{-\Theta^{1-\beta}/12e} = \exp \Bigl(-\frac{\Theta^{1-\beta}}{12e}\log X\Bigr)
}
and for $r=0$, 
\eqs{T_0^{-\ell_0} < 4^{-\ell_0/12} < 4^{- (\log_2 X)^{\alpha\beta}/12} = \exp \Bigl(- (\log_2 X)^{\alpha\beta} \frac{\log 4}{12} \Bigr).
}
Since the largest error is introduced when $r=0$, we see that the total error of removing the restrictions for all $0 \le r \le J$ is $\ll J T^{-\ell_0} = o (X\log X)$, since $J < \log_3 X$ and $\alpha\beta> 2$. 

Finally, we write the sums over $a_r$ without the restrictions as
\bs{
	&\sum_{c_r \in \cI_r}  \frac 1  {\N(c_r)^{3/2}}  \pr{\pi \equiv 1 \mod 3\\ \pi \mid c_r}  \frac  {\N(\pi) } {\N(\pi) + 2} 
	\su{a_r, b_r \in \cI_r\\ a_rb_r = c_r^3} \lambda(a_r)f(a_r,J)\nu(a_r) \\
	&= \prod_{\xi \in \cI_r} \Bigl(1 + \frac  {\N(\xi) } {\N(\xi) + 2}  \sum_{n \ge 1}   \frac 1  {\N(\xi )^{3n/2}} 
	\su{0 \le m \le 3n} \lambda(\xi^m)f(\xi^m,J)\nu(\xi^m) \Bigr)\\
	&= \prod_{\xi \in \cI_r} \Bigl(1 + \frac  {\N(\xi) } {\N(\xi) + 2}  \su{m \ge 0} \frac{(-f(\xi,J))^m}{m!} \su{n \ge \max(1, m/3)}   \frac 1  {\N(\xi )^{3n/2}} \Bigr)\\
	&= \prod_{\xi \in \cI_r} \Bigl(1 + \frac  {\N(\xi) } {\N(\xi) + 2}  \frac 1 {(1-\N(\xi)^{-3/2})} \sum_{m \ge 0} (-1)^m b_m \Bigr),
}
where
\eqs{b_m = \frac{f(\xi, J)^m}{m!\N(\xi)^{\frac 3 2 \max(1,\ceil{m/3})}} \ge 0,}
with equality when there is a prime $\xi$ with $\N(\xi) = X^\theta_J$. This case gives the same result as the sum over $b$. Hence, we shall assume that $\N(\xi) < X^{\theta_J}$. Since for each $m \ge 0$, 
\eqs{
	\frac{b_{m+1}}{b_m} = \frac{f(\xi, J) \N(\xi)^{3\max(1,\ceil{m/3})/2}} {(m+1)\N(\xi)^{3\max(1,\ceil{(m+1)/3})/2}} < 1,
}
and that $b_m < \frac 1 {m!}$, we have a convergent alternating series. Therefore, 
\eqn{ \label{maintermlowerbound2}
b_0 = \frac 1 {\N(\xi)^{3/2}} >  \sum_{m \ge 0} (-1)^m b_m > b_0 - b_1 = \frac {1 - f(\xi,J) } {(\N(\xi))^{3/2}} > 0.
}

Combining \eqref{maintermlowerbound1} and \eqref{maintermlowerbound2}, we see that the main term of the first mollified moment is $C_X X\log X$, where
\multn{\label{CX}
C_X = c_0 \pr{\xi \equiv 1 \mod 3\\ \N\xi \ge X^{\theta_J}} \biggl(1 + \frac {\N(\xi)}{(\N(\xi) + 2)(\N(\xi)^{3/2}-1)} \biggr) \\
\cdot \pr{ \xi \equiv 1 \mod 3\\ k_0 < \N\xi  < X^{\theta_J} } \Bigl(1 + \frac {\N(\xi)}{(\N(\xi) + 2)(\N(\xi)^{3/2}-1)} \sum_{m \ge 0} (-1)^m b_m \Bigr)
}
and that 
\eqs{\label{CXlowerbound}
	C_X  > c_0 = \frac {4\pi^2} {2187(1-1/\sqrt 3)}  \pr{\pi \equiv 1 \mod 3}  \Bigl( 1 - \frac 3 {\N(\pi)^2} + \frac 2 {\N(\pi)^3}
	\Bigr),
}
and
\eqn{ \label{c1}
	C_X  < c_1 := c_0  \pr{\xi \equiv 1 \mod 3\\ \N\xi \ge X^{\theta_J}} \biggl(1 + \frac {\N(\xi)}{(\N(\xi) + 2)(\N(\xi)^{3/2}-1)} \biggr).
}
This proves the claimed result in Theorem \ref{thm:firstmoment}, provided that $X$ is sufficiently large in terms of $\alpha, \beta, \eps$ and \eqref{ThetaCond1stmoment} is satisfied.

\subsection{Proof of Theorem \ref{thm:posprop}}
We take 
\eqs{k= 2, \quad \kappa = 1, \quad a= 1.0299, \quad \beta = 0.916, \quad \epsilon = 10^{-5}, \quad \alpha = 7, }
and 
\eqs{\Theta = 5.8 0 2 5 9 3 5 5 1 5 \times  10^{-44}.}
Then, the conditions in \eqref{alphabetacond}, \eqref{eta0bound}, \eqref{etaJcond}, \eqref{thetacond1}, \eqref{etajcond}, \eqref{betacond}, \eqref{ThetaBetacond}, \eqref{thetacond4}, \eqref{thetacond5} and \eqref{ThetaCond1stmoment}  are all satisfied. Furthermore, we have 
\eqs{R_1 \approxeq  -4.7107876828 \times 10^{40},\qquad R_2 \approxeq 2.8043085602 \times 10^{42},}
which were defined in \eqref{R1R2}. Since $e^r \ge r$, we have for $r \ge 0$ that 
\eqs{f(r) = e^r (R_1 - R_2 r) + 8 r \le r( R_1 - R_2r + 8	) \le  (R_1+8)r.}
Thus, we see that
\bs{
\su{0 \le r \le J-1} (r+1) &\exp \biggl( f(r) + O \Bigl( \frac {(\log_2 X)^{2\alpha}} {\log X} \Bigr)  \biggr) \\
&< \sum_{r \ge 0 } (r+1)e^{r (R_1+8)} = \frac 1 {(1-e^{R_1 + 8})^2}
}
for sufficiently large $X$. Moreover, with our choices, 
\eqs{\cD_2 \le 2.6176409874 \times 10^{15}.}
This leads to the inequality 
\bs{
\eqref{HMsplit1} & <  |\sF(X)| ( 1 + o_X(1)) 2.6176409874 \times 10^{15} \( 2 \exp \bigl( 2e^{1/4} + 2e^{1/8} \bigr)  	
+ \exp \bigl(2e/\Theta \bigr) \) \\
&= |\sF(X)| ( 1 + o_X(1)) \Bigl(6.5837089699\times 10^{17} + \exp \bigl(\exp ( 101.248586291 )\bigr)\Bigr).
}
We conclude that the estimate
\eqs{
\sum_{c \in \sF(X) \cap \cT_0} |L(1/2,\chi_c)|^2 |M(c, 1)|^2 < e^{e^{101.3}} |\sF(X)| (1 + o_X(1))
}
holds for sufficiently large $X$. Furthermore, it follows from Lemma \ref{familysize} with $n = 1$ that 
\eqs{|\sF(X)| \le (|C_1(1)| + o_X(1))X \log X, \quad |C_1(1)| \le \frac{4\pi^2 g_{\psi_0} (1) }{2187} = (1 - 1/\sqrt 3)c_0.}
Combining this result with \eqref{CXlowerbound} and using Cauchy-Schwarz inequality it follows that 
\eqs{
\frac 1 {|\sF(X)|}\su{c\in \sF(X) \\ L(1/2, \chi_c) \neq 0} 1 >
\frac{c_0^2 + O(1/\log X)}{e^{e^{101.3}}(1+o_X(1)) ((1 - 1/\sqrt 3)c_0 + o_X(1))^2}.
}
The claim of Theorem \ref{thm:posprop} follows upon letting $X$ tend to infinity.

\nocite{*}

\bibliographystyle{amsplain}

\end{document}